\documentclass[preprint,12pt]{elsarticle}
\usepackage{amsfonts}

\usepackage{amssymb,color}

\journal{AMM}

\begin{document}

\begin{frontmatter}
\newtheorem{theo}{Theorem}
\newtheorem{lemma}[theo]{Lemma}
\newtheorem{remark}[theo]{Remark}
\newtheorem{cor}[theo]{Corollary}
\newtheorem{proof}[theo]{Proof}
\newtheorem{definition}[theo]{Definition}


\title{Interpolated variational iteration method for initial value problems}


\author{Davod Khojasteh Salkuyeh$^\dag$\footnote{Corresponding author} and Ali Tavakoli$^\S$ }

\address{$^\dag$Faculty of Mathematical Sciences, University of Guilan, Rasht, Iran\\
e-mail: khojasteh@guilan.ac.ir, salkuyeh@gmail.com\\[2mm]
$^\S$Mathematics Department, Vali-e-Asr University of Rafsanjan, Iran.\\
e-mail: tavakoli@mail.vru.ac.ir}

\begin{abstract}
\indent  In order to solve an initial value problem by the variational iteration method, a sequence of functions is produced which  converges to the solution under some suitable conditions. In the nonlinear case, after a few iterations the terms of the sequence become complicated,  and therefore, computing a highly accurate solution would be difficult
or even impossible. In this paper, for one-dimensional initial value problems, we propose a new
approach which is based on approximating each term of the sequence by a piecewise linear function.
Moreover, the convergence of the method is proved.  Three illustrative examples are given to show the superiority
of the proposed method over the classical variational iteration method.
\end{abstract}


\end{frontmatter}


\section{Introduction} \label{SEC1}

The He's variational iteration method (VIM) \cite{He1,He4} is a powerful mathematical technique to solve linear and nonlinear problems that can be implemented easily in practice. It has been successfully applied for solving various ODEs and PDEs \cite{He1,HeIJMPB,He4,He7,He8,He9}. Convergence of the VIM has been investigated in some papers \cite{SalkuyehCAMWA,Tatari}.  In \cite{Zhao}, Zhao and Xiao applied the VIM for solving singular perturbation initial value problems and investigated its convergence.  Recently, Wu and Baleanu in \cite{WuBaleanu2} have introduced a new way  to define the Lagrange multipliers in the VIM for solving  fractional differential equations with the Caputo derivatives. They also developed the VIM to $q$-fractional difference equations in \cite{WuBaleanu1}. The review article \cite{GCWu} concerns some new applications of the VIM to numerical simulations  of differential equations and fractional differential equations.

Consider the following differential equation
\begin{equation}\label{Eq001}
 \mathcal{L}u(t)+\mathcal{N}u(t)=g(t),
\end{equation}
\noindent where $\mathcal{L}$ and $\mathcal{N}$ are, respectively,  linear and nonlinear operators and  $g$ is an inhomogeneous term.
Given an initial guess $u_0(t)$, the VIM to solve (\ref{Eq001}) takes the following form
\begin{equation}\label{Eq002}
 u_{m+1}(t)= u_m(t)+\int_{t_0}^{t} \lambda(s,t) (\mathcal{L}u_n(s)+
\mathcal{N}u_m(s)-g(s))ds,\quad m=0,1,2,\ldots.
\end{equation}
where  $\lambda$ is a general Lagrange multiplier (for more details see \cite{He1,He4,He7}).

 As we see the VIM generates a sequence of functions that can converge to the solution of the problem under some suitable conditions. After a few iterations each term of this sequence often involves a definite integral whose integrand contains several nonlinear terms. In this case, computing a high accuracy solution would be difficult or even impossible by the common softwares such as MAPLE, MATHEMATICA or MATLAB.  See \cite{Abbasbandy} for an  application of the VIM to quadratic Riccati differential equations that illustrates this. To overcome  this problem, Geng et al. in \cite{Geng} have introduced a piecewise variational iteration method for the quadratic Riccati differential equation. Their main idea is that for solving the problem in a large interval, the interval is split into several subintervals and the  problem is solved in each subinterval by the VIM in a progressive manner. In this way, the accuracy of the computed solution can be considerably improved. However, the main mentioned problem is still unsolved, since in each subinterval the VIM is directly implemented.

In this paper, to solve one-dimensional initial value problems, we modify the VIM in such a way that each term of the sequence is interpolated  by a piecewise linear function. Hereafter, the proposed method is called the IVIM (for interpolated VIM). In spite of the VIM, the IVIM does not need any symbolic computation and all of the computations are done numerically. Therefore, we can compute hundreds of the sequence terms in a small amount of time.

The rest of the paper is organized as follows.  The IVIM is introduced in Section \ref{SEC3}. Convergence of the proposed method is discussed in Section \ref{SEC3}. Section \ref{SEC4} is devoted to presenting three numerical examples. Some concluding remarks are given in Section \ref{SEC5}.

\section{The IVIM }  \label{SEC2}

Consider the one-dimensional initial value problem
\begin{equation}\label{Eq003}
\left\{
  \begin{array}{ll}
    u'(t)=f(t,u(t)), & t\in [a,T], \\
    u(a)=u_a.
  \end{array}
\right.
\end{equation}
For simplicity, we assume that $u_a=0$, otherwise one can use a simple change of variable $\tilde{u}=u-u_a$ to obtain $\tilde{u}(a)=0$. In this case, Eq. (\ref{Eq002}) becomes
\begin{equation}\label{Eq004}
 u_{m+1}(t)= u_m(t)+\int_{a}^{t} \lambda(s,t) \left( u_m'(s)-f(s,u_m(s))\right)ds,\quad m=0,1,2,\ldots,
\end{equation}
where $u_0(t)$ satisfies the initial condition $u_0(a)=0$. By integration by parts, Eq. (\ref{Eq004}) can be rewritten as
\begin{equation}\label{Eq005}
 u_{m+1}(t)=  G_m(t)-\int_{a}^{t}H_m(s,t)ds,
\end{equation}
where
\begin{eqnarray*}
 G_m(t)   \hspace{-0.2cm}&=&\hspace{-0.2cm} (1+\lambda(t,t))u_m(t)-\lambda(a,t)u_m(a),
  \end{eqnarray*}
 and
 \begin{eqnarray*}
  H_m(s,t) \hspace{-0.2cm}&=&\hspace{-0.2cm} \frac{\partial \lambda}{\partial s}(s,t)~u_m(s)+\lambda(s,t)f(s,u_m(s)).
 \end{eqnarray*}

In order to present the IVIM, we take a natural number $n$ and  discretize the interval $[a,T]$ into  $n-1$ subintervals with step size $h=(T-a)/(n-1)$ and grid points
\[
t_i=a+(i-1)h,\quad i=1,2,\ldots,n.
\]
Now, we define the well-known B-spline basis functions (see \cite{Cheney})  of first-order on the nodal points $t_i$, i.e.
\begin{eqnarray*}
 \varphi_i(t)  &=&  \left\{
                  \begin{array}{ll}
                  \displaystyle\frac{t-t_{i-1}}{h}, & t_{i-1} \leq t < t_i, \\[4mm]
                  \displaystyle\frac{t_{i+1}-t}{h}, & t_{i} \leq t \leq t_{i+1}, \hspace{2cm} i=2,3,\ldots,n-1,\\[4mm]
                  \quad 0 , &  t<t_{i-1}~\textrm{or}~t>t_{i+1},\\[4mm]
                  \end{array}
                  \right. \\
   \varphi_n(t) &=&\left\{
                   \begin{array}{ll}
                     \displaystyle\frac{t-t_{n-1}}{h}, & t_{n-1} \leq t \leq t_n, \\[4mm]
                     \quad 0 , & \quad t<t_{n-1},
                   \end{array}
                  \right. \\
\end{eqnarray*}
on $[a,T]$. Each function $\varphi_i$ is equal to $1$ at grid point $t_i$, and equal to zero at other grid points. Let
\[
X_h=\textrm{span}\{\varphi_{i}:i=2,3,\ldots,n\}.
\]
Every $v^{(h)}$ in $X_h$ is a piecewise linear function of the form
\[
v^{(h)}(t)=\sum_{i=2}^n \alpha_i \varphi_i (t),\quad t\in [a,T].
\]
Obviously, $\alpha_i$ is the value of  $v^{(h)}$ at the grid point $x_i$, i.e., $\alpha_i=v^{(h)}(t_i)$. Furthermore, $v^{(h)}(a)=0$, since $\varphi_i (a)=0$, $i=2,3,\ldots,n$. See \cite{Reddy} for more details about the properties of this kind of functions.

We now return to Eq. (\ref{Eq005}).  The aim is to recursively compute  $u_{i}$'s. Let $u_m$ has been computed and we intend to compute $u_{m+1}$.
The main idea of our method is to compute a piecewise linear interpolation $u_{m+1}^{(h)}$  of $u_{m+1}$  in $X_h$, instead of computing $u_{m+1}(t)$.
In fact, instead of $u_{m+1}$, we consider its  projection onto $X_h$. To do so, we have
\begin{equation}\label{Eq006}
u_{m+1}^{(h)}(t)=\displaystyle \sum_{r=2}^n u_{m+1}(t_r)\varphi_{r}(t) ,\quad t\in [a,T].
\end{equation}
All we need here is to compute $u_{m+1}(t_i)$, $i=2,3,\ldots,n$. By using Eq. (\ref{Eq005}), we have
\begin{equation}\label{Eq007}
u_{m+1}(t_i)= G_m(t_i)-\int_{a}^{t_i}H_m(s,t_i)ds.
\end{equation}
Now, we can approximate the integral term in the latter equation by different methods of numerical integration. But, we  again compute a piecewise linear interpolation of $H_m(s,t_i)$ in $X_h$ as
\[
H_m^{(h)}(s,t_i)=\sum_{r=2}^n H_m(t_r,t_i)\varphi_r(s),\quad s\in [a,t_i].
\]
Replacing $H_m(s,t_i)$ by $H_m^{(h)}(s,t_i)$ in Eq. (\ref{Eq004}) yields
\begin{equation}\label{Eq008}
u_{m+1}(t_i)  \approx  G_m(t_i)-\sum_{r=2}^n H_m(t_r,t_i)\int_{a}^{t_i} \varphi_r(t)dt =G_m(t_i)-\sum_{r=2}^n H_m(t_r,t_i)\mu_{ri},
\end{equation}
where $\mu_{ri}=\displaystyle\int_{a}^{t_i} \varphi_r(t)dt$,  $r=2,3,\ldots,n$. It is easy to see that
\[
\mu_{ri}=\left\{
           \begin{array}{ll}
            0,           &  i\leq r-1, \\[2mm]
           \displaystyle \frac{h}{2},  &  i=r,  \\ [4mm]
           \displaystyle h , & i\geq r+1.
           \end{array}
         \right.
\]
Now, from this and Eq. (\ref{Eq008}) we get
\begin{equation}\label{Eq009}
u_{m+1}(t_i)\approx \hat{u}_{m+1}(t_i)=G_m(t_i)-h\sum_{r=2}^{i-1} H_{m}(t_r,t_i)-\frac{h}{2}H_m(t_i,t_i),\quad i=2,3,\ldots,n.
\end{equation}
Then, from (\ref{Eq006}) we set
\begin{equation}\label{Eq010}
u_{m+1}^{(h)}(t)\approx \displaystyle \sum_{r=2}^n \hat{u}_{m+1}(t_r)\varphi_{r}(t).
\end{equation}

\noindent {\bf Remark 1.} Consider the special case that $\mathcal{L}u\equiv u'+\alpha u$, for some $\alpha \in \mathbb{R}$. This means that  $f(t,u(t))=-\alpha u(t)-\mathcal{N}u(t)+g(t)$, where $\mathcal{N}u(t)$ is the nonlinear term appeared in  Eq. (\ref{Eq001}).  In this case, we have $\lambda(s,t)=-e^{\alpha(s-t)}$, and therefore $1+\lambda(t,t)=0$. Obviously, from Eq. (\ref{Eq004}) we see that $u_m(a)=0$, for $m=0,1,\ldots$. Hence, $G_m(t)=0$, for $m=0,1,\ldots$. On the other hand, we have
\begin{eqnarray}
\nonumber H_m(s,t) \hspace{-0.2cm}&=&\hspace{-0.2cm} e^{\alpha(s-t)}(-\alpha u_m(s)+f(s,u_m(s))) \\
\nonumber         \hspace{-0.2cm}&=&\hspace{-0.2cm} e^{\alpha(s-t)}(-\alpha u_m(s)+\alpha u_m(s)+\mathcal{N}u_m(s)-g(s))\\
         \hspace{-0.2cm}&=&\hspace{-0.2cm} e^{\alpha(s-t)}(\mathcal{N}u_m(s)-g(s)). \label{Eq011}
\end{eqnarray}

\section{Convergence of the method}\label{SEC3}

In this section, we present the convergence of the VIM and the IVIM. We first recall a lemma which gives conditions for the existence and uniqueness of a solution to the problem (\ref{Eq003}).

\begin{lemma}
Consider the initial-value problem (\ref{Eq003}), where $f:[a,T]\times \Bbb R\rightarrow\Bbb R$ is continuous with respect to the first variable and satisfies a Lipschitz condition with respect to the second variable, i.e.,
\begin{equation}\label{conv00}
|f(t,u(t))-f(t,v(t))|\leq \gamma(t) |u(t)-v(t)|,\quad t\in [a,T],
\end{equation}
where $\gamma(t)$ is a nonnegative continuous function. Then, there exists a unique solution to the problem.
\end{lemma}

\noindent{\bf Proof}. See  \cite{buc,tav}.\\

To establish convergence of  the VIM, we state the following theorem.

\begin{theo}\label{theo1}
Let $f:[a,T]\times C^1([a,T])\rightarrow\Bbb R$ be a  continuous function with respect to its first variable which satisfies the Lipschitz condition (\ref{conv00}).
If $u_m(t)\in C^1([a,T])$ for $m=0,1,2,\ldots$, then the sequence defined by (\ref{Eq004}) converges to the solution of (\ref{Eq003}).
\end{theo}

\noindent {\bf Proof}.
By Eq. (\ref{Eq004}), we have
\begin{equation}\label{conv01}
 u_{m+1}(t)= \int_{a}^{t} ((1+\lambda(s,t))  u_m'(s)-\lambda(s,t)f(s,u_m(s)))ds,\quad m=0,1,2,\ldots,
\end{equation}
since $u_m(a)=0$.
Let $u(t)$ be the exact solution of Eq. (\ref{Eq003}). By Eq. (\ref{conv01}), it is readily seen that
\begin{eqnarray}\label{conv02}
 u(t)= \int_{a}^{t} ((1+\lambda(s,t))  u'(s)-\lambda(s,t)f(s,u(s)))ds.
 \end{eqnarray}
Let $E_m(t)=u_m(t)-u(t)$ be the error produced by the VIM at $m$th iteration.
 Since the first derivative of $E_0(t)$ is continuous on the closed interval $[a,T]$, there  exists a constant $M$ such that
 \begin{eqnarray}\label{conv04}
 |E_0'(t)|\leq M|E_0(t)|,\qquad \forall t\in[a,T].
  \end{eqnarray}
From  Eqs. (\ref{conv01}), (\ref{conv02}) and (\ref{conv00}), it follows that
\begin{eqnarray}
\nonumber |E_{1}(t)| \hspace{-0.2cm}&=&\hspace{-0.2cm} |\int_{a}^{t} ((1+\lambda(s,t))  E_0'(s)-\lambda(s,t)(f(s,u_0(s))-f(s,u(s))))ds|  \\
       \hspace{-0.2cm}& \leq &\hspace{-0.2cm}  \int_{a}^{t} |1+\lambda(s,t)||E_0'(s)| ds+\int_{a}^{t}|\gamma(s)\lambda(s,t)||E_0(s)| ds.  \label{conv05}
\end{eqnarray}
Then, it turns out that by Eqs. (\ref{conv04}) and (\ref{conv05}),
\begin{equation}\label{conv06}
|E_{1}(t)|\leq M_a \int_{a}^{t} |E_0(s) | ds,
\end{equation}
in which
\[
M_a=\max\large\{M \max_{s,t\in[a,T]} |1+\lambda(s,t)|, \max_{s,t \in[a,T]}|\gamma(s)\lambda(s,t)| \large\}.
\]
Let $\|E_i\|=\displaystyle\max_{t\in [a,T]} |E_i(t)|,$ $ i=0,1,\ldots$.  Hence,
\begin{equation}\label{conv07}
|E_{1}(t)|\leq  \|E_{0}\| M_a(t-a).
\end{equation}
Therefore, similar to Eq. (\ref{conv06}) and using Eq. (\ref{conv07}), we deduce
\begin{eqnarray*}
       |E_{2}(t)|    \hspace{-0.2cm}& \leq &\hspace{-0.2cm} M_a\int_{a}^{t}|E_1(s)| ds \leq M_a^2 \|E_{0}\|\int_{a}^{t}(s-a) ds  =  \|E_{0}\|\frac{M_a^2(t-a)^2}{2!},\\
       |E_{3}(t)|    \hspace{-0.2cm}& \leq &\hspace{-0.2cm} M_a\int_{a}^{t}|E_2(s)| ds \leq M_a^3 \|E_{0}\|\int_{a}^{t}(s-a)^2 ds  =  \|E_{0}\|\frac{M_a^3(t-a)^3}{3!},
\end{eqnarray*}
and in general,
\begin{eqnarray*}
|E_{m}(t)|\leq \|E_{0}\| \frac{M_a^m(t-a)^m}{m!}\leq \|E_{0}\|\frac{M_a^m(T-a)^m}{m!}.
\end{eqnarray*}
Then $\|E_{m}\|\rightarrow 0$ as $m$ tends to $\infty$  and the convergence of the method is established.\qquad $\Box$

\bigskip

Now, we are ready to prove the convergence of the IVIM.

\bigskip

\begin{theo}
 Let $f:[a,T]\times C^1([a,T])\rightarrow\Bbb R$ be two times  continuously differentiable function with respect to the first variable and satisfies the Lipschitz condition (\ref{conv00}). If the VIM of (\ref{Eq004}) is convergent, then  so is the IVIM of (\ref{Eq010}).
\end{theo}
{\bf Proof.}
Define $\bar{u}_{m}(t):=\sum\limits_{i=2}^n \hat{u}_{m}(t_i)\varphi_{i}(t)$. For each grid  point $t_r$,  we have $\bar{u}_{m}(t_r)=\hat{u}_{m}(t_r)$.  Also, by Eq. (\ref{Eq009}), we have
\begin{eqnarray}
  \nonumber        \hat{u}_{m+1}(t_i) \hspace{-0.2cm}&= &\hspace{-0.2cm} (1+\lambda(t_i,t_i))u_m(t_i)-h\sum\limits_{r=2}^{i-1} (\frac{\partial\lambda}{\partial s}(s,t_i)\mid_{s=t_r}u_m(t_r)+\lambda(t_r,t_i)f(t_r,u_m(t_r)))\\
                              \hspace{-0.2cm}&  &\hspace{-0.2cm} -\frac{h}{2}((\frac{\partial\lambda}{\partial s}(s,t_i)\mid_{s=t_i}u_m(t_i)+\lambda(t_i,t_i)f(t_i,u_m(t_i))). \label{conv08}
\end{eqnarray}
On the other hand, by Eq. (\ref{Eq005}) it is clear that the exact solution $u(t)$ at $t=t_i$ satisfies the following relation
\begin{eqnarray}\label{conv09}
u(t_i)=(1+\lambda(t_i,t_i))u(t_i)-\int\limits_a^{t_i}(\frac{\partial \lambda}{\partial s}(s,t_i)~u(s)+\lambda(s,t_i)f(s,u(s)))ds.
\end{eqnarray}
 Now, by the trapezoidal rule of integration, we have\\[4mm]
$\displaystyle \int_a^{t_i}(\frac{\partial \lambda}{\partial s}(s,t_i)~u(s)+\lambda(s,t_i)f(s,u(s)))ds$
\begin{eqnarray}
\nonumber \qquad \hspace{-0.2cm}& = &\hspace{-0.2cm} \frac{h}{2}(\frac{\partial\lambda}{\partial s}(s,t_i)\mid_{s=t_i}u(t_i)+\lambda(t_i,t_i)f(t_i,u(t_i))+\lambda(a,t_i)f(a,0)) \\
\nonumber \hspace{-0.2cm}&  &\hspace{-0.2cm} +h\sum_{r=2}^{i-1} (\frac{\partial\lambda}{\partial s}(s,t_i)\mid_{s=t_r}u(t_r) +\lambda(t_r,t_i)f(t_r,u(t_r)))\\
 \hspace{-0.2cm}&  &\hspace{-0.2cm}  -\frac{t_i-a}{12}h^2 \frac{d^2}{ds^2}f(s,u(s))|_{s=\xi}, \label{conv10}
\end{eqnarray}
where $\xi\in(a,t_i)$. Let $\bar{E}_m(t_i)=\bar{u}_{m}(t_i)-u(t_i)=\hat{u}_{m}(t_i)-u(t_i)$ and $E_m(t_i)=u_m(t_i)-u(t_i)$. By Eqs. (\ref{conv08}), (\ref{conv09}) and (\ref{conv10}), we get
\begin{eqnarray*}
\nonumber |\bar{E}_{m+1}(t_i)| \hspace{-0.2cm}& = &\hspace{-0.2cm} |\hat{u}_{m+1}(t_i)-u(t_i)|  \\
\nonumber \hspace{-0.2cm}& =  &\hspace{-0.2cm}   |(1+\lambda(t_i,t_i))E_m(t_i)-h\sum\limits_{r=2}^{i-1} (\frac{\partial\lambda}{\partial s}(s,t_i)\mid_{s=t_r}E_m(t_r)\vspace{.3cm}\\
\nonumber \hspace{-0.2cm}&   & \hspace{-0.2cm} +\lambda(t_r,t_i)(f(t_r,u_m(t_r))-f(t_r,u(t_r))))
-\frac{h}{2}(\frac{\partial\lambda}{\partial s}(s,t_i)\mid_{s=t_i}E_m(t_i)\\
\nonumber \hspace{-0.2cm}&   &\hspace{-0.2cm}+\lambda(t_i,t_i)(f(t_i,u_m(t_i))-f(t_i,u(t_i))))-\frac{h}{2}\lambda(a,t_i)f(a,0)\\
\hspace{-0.2cm}&   &  \hspace{-0.2cm} +\frac{t_i-a}{12}h^2 \frac{d^2}{ds^2}f(s,u(s))|_{s=\xi}|.
\end{eqnarray*}
Now, by the Lipschitz condition, it follows that
\begin{eqnarray*}
\nonumber|\bar{E}_{m+1}(t_i)| \hspace{-0.2cm}& \leq &\hspace{-0.2cm} |1+\lambda(t_i,t_i)|~|E_m(t_i)|+h\sum_{r=2}^{i-1} (|\frac{\partial\lambda}{\partial s}(s,t_i)\mid_{s=t_r}|~ |E_m(t_r)|\\
\nonumber\hspace{-0.2cm}&&\hspace{-0.2cm} +|\gamma(t_r)\lambda(t_r,t_i)|~|E_m(t_r)|)
+\frac{h}{2}(|\frac{\partial\lambda}{\partial s}(s,t_i)\mid_{s=t_i}|~|E_m(t_i)|\\
\nonumber\hspace{-0.2cm}&&\hspace{-0.2cm}+|\gamma(t_i)\lambda(t_i,t_i)|~|E_m(t_i)|+\frac{h}{2}|\lambda(a,t_i)f(a,0)|
+\frac{t_i-a}{12}h^2 | \frac{d^2}{ds^2}f(s,u(s))|_{s=\xi}|\\
\nonumber\hspace{-0.2cm}&\leq&\hspace{-0.2cm}  K_1|E_m(t_i)|+hK_2\sum_{r=2}^{i-1} |E_m(t_r)|+\frac{h}{2}K_2 |E_m(t_i)|\\
\hspace{-0.2cm}&&\hspace{-0.2cm} +
\frac{h}{2}|\lambda(a,t_i)f(a,0)|+\frac{t_i-a}{12}h^2K_3,
\end{eqnarray*}
where $|1+\lambda(t,t)|\leq K_1$, $|\frac{\partial\lambda}{\partial s}(s,t)|+|\gamma(t)\lambda(s,t)|\leq K_2$
and $|\frac{d^2}{ds^2}f(s,u(s))|\leq K_3$ for some constants $K_1,K_2$ and $K_3$. Now, if the VIM is convergent, then $|E_m(t_i)|\rightarrow 0$  as $m\rightarrow\infty$ and  therefore, as a result when $m\rightarrow \infty$ and $h\rightarrow 0$, we have $|\bar{E}_{m}(t_i)|\rightarrow 0$ for each grid point $t_i$.
It is necessary to mention that
\[
\lim_{m\rightarrow\infty}\lim_{h\rightarrow 0} h\sum_{r=2}^{i-1}|E_{m}(t_r)|=\lim_{m\rightarrow\infty}\int_{x_1}^{x_{i-1}} |E_{m}(t)|dt=0.
\]
Therefore the proof is completed.
 \qquad $\Box$

\section{Illustrative examples}\label{SEC4}

In this section, we present three numerical examples to show the effectiveness of our method.  In order to do so, we compare the results of the IVIM to those of the VIM.

\bigskip

\begin{figure}
\centering
\includegraphics[height=6cm,width=8cm]{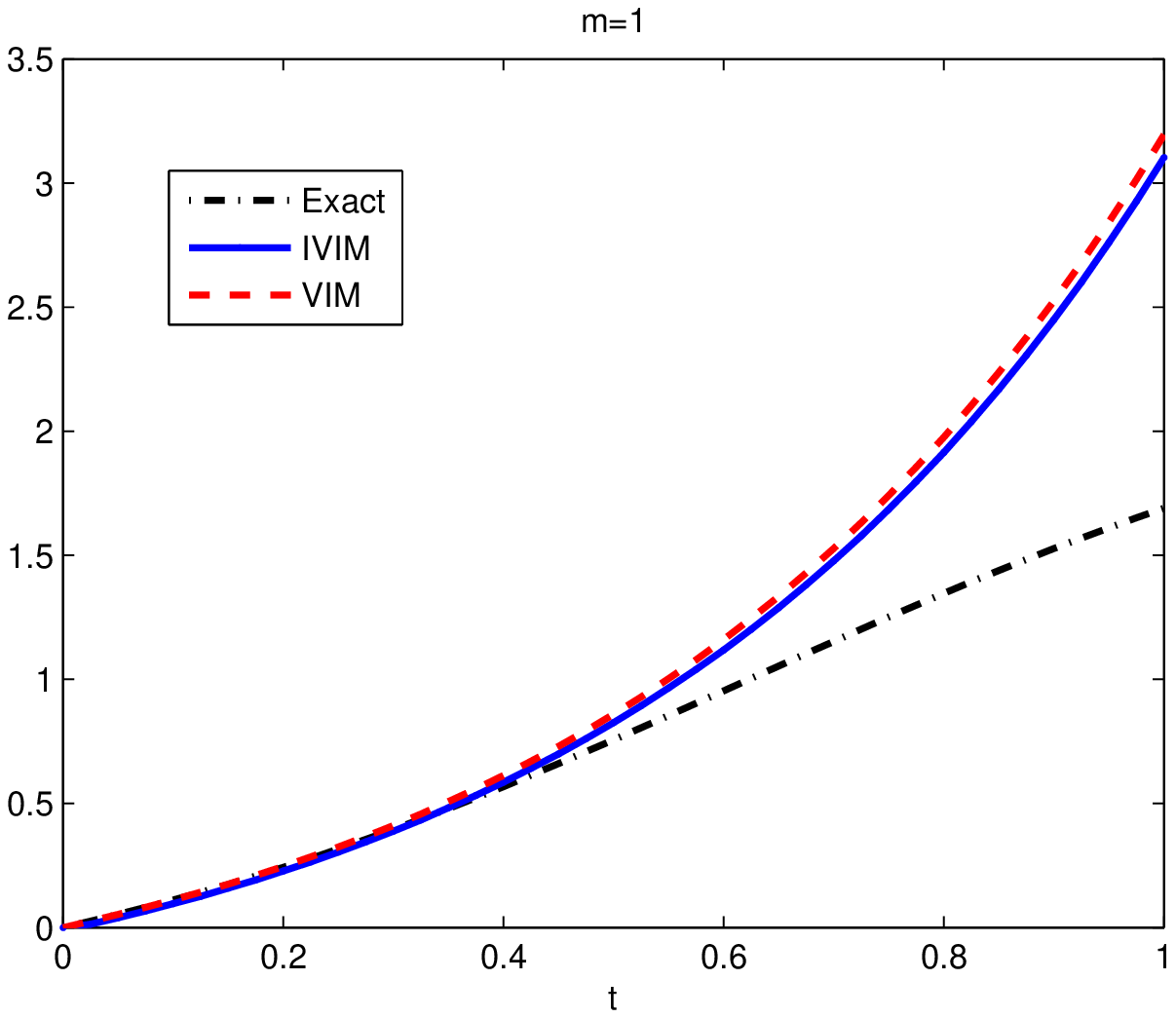}\includegraphics[height=6cm,width=8cm]{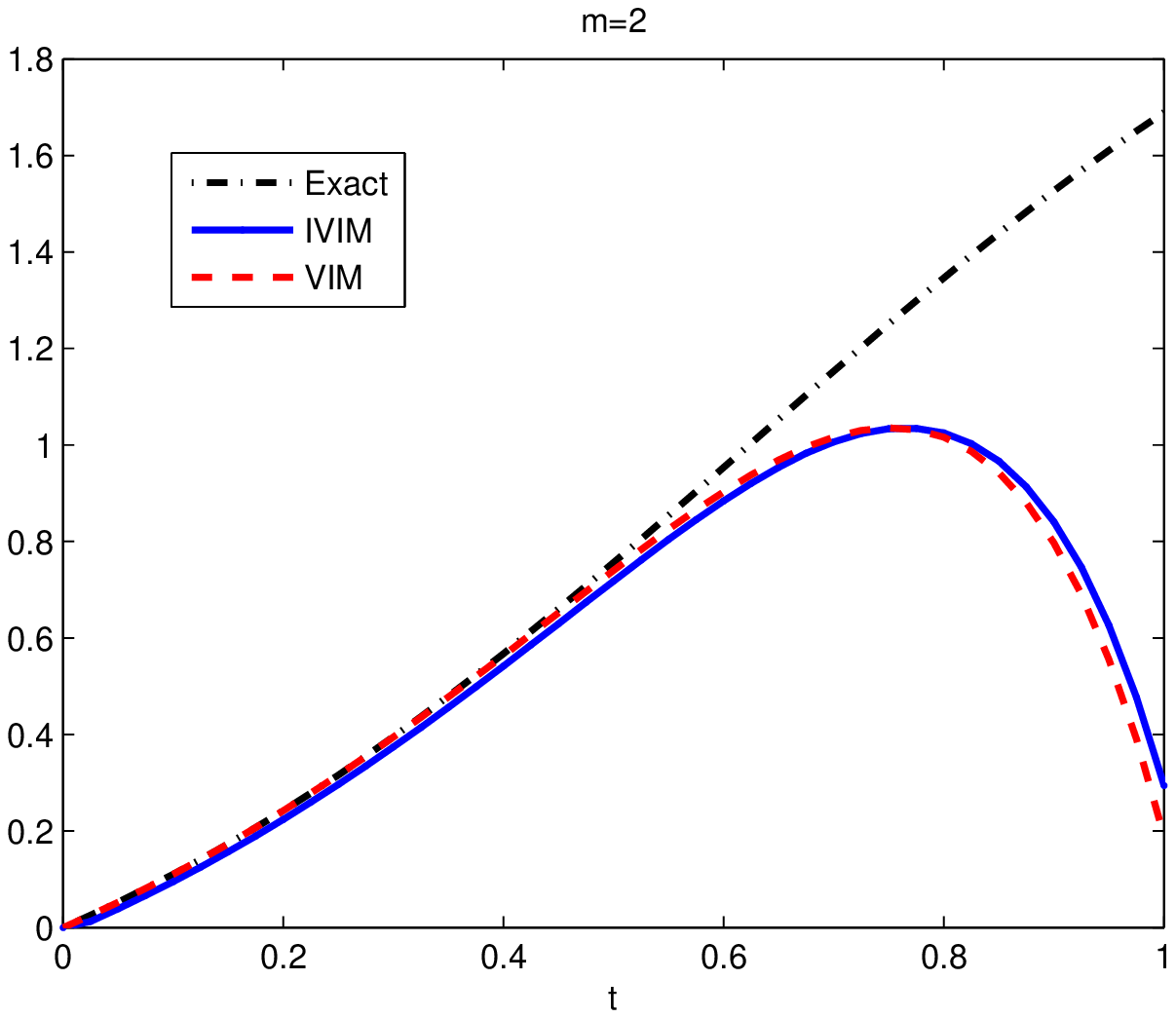}
\includegraphics[height=6cm,width=8cm]{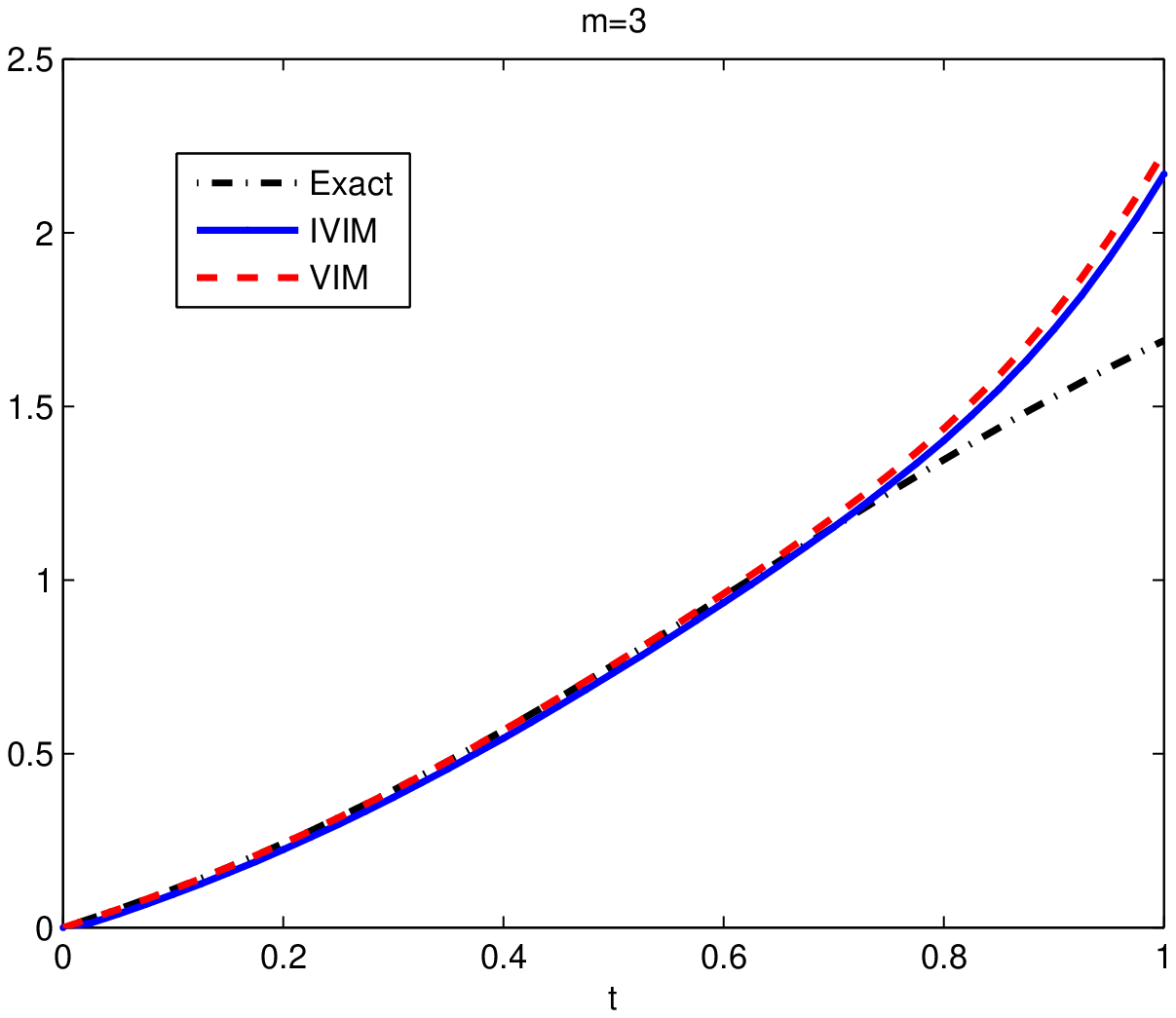}\includegraphics[height=6cm,width=8cm]{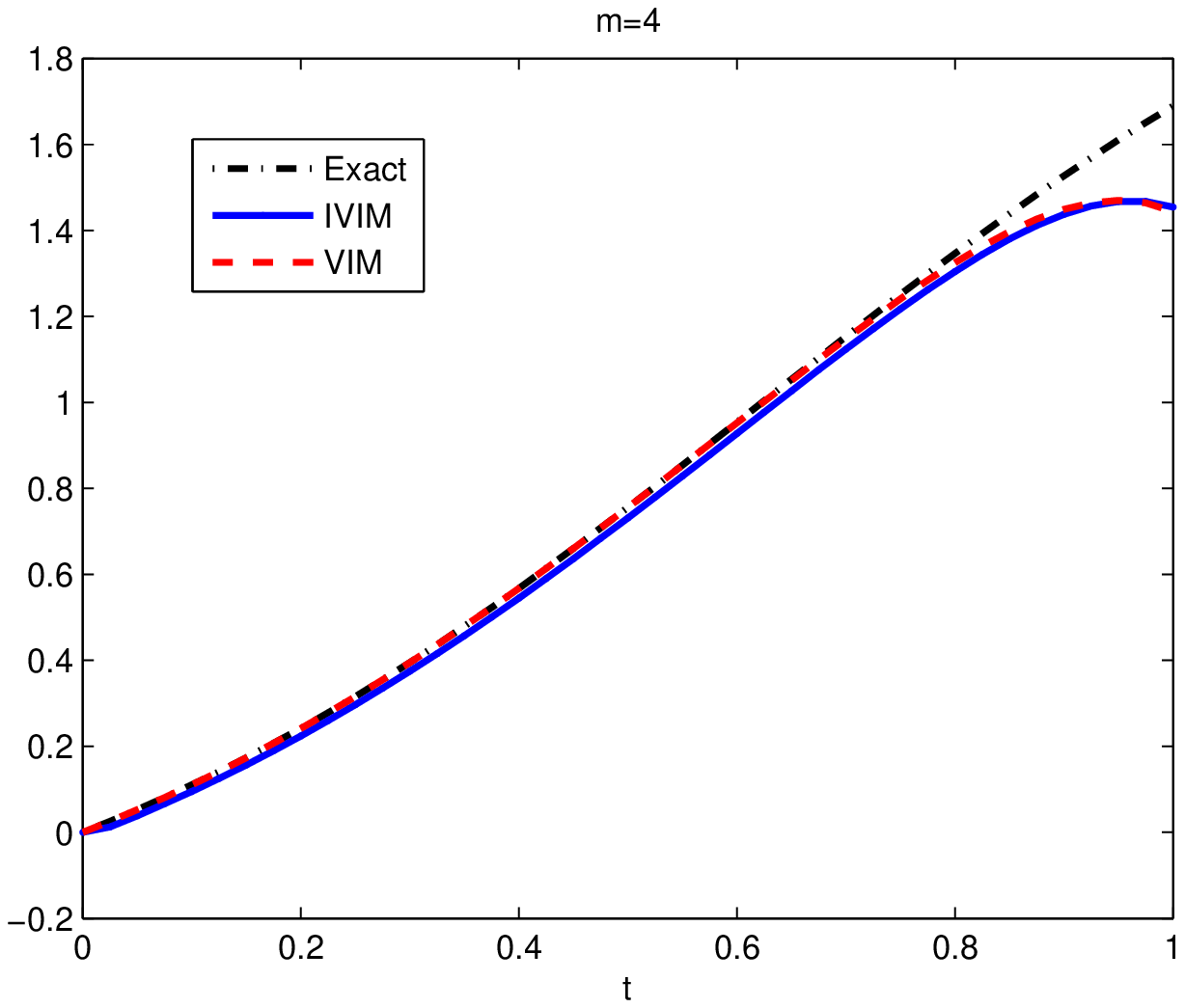}
\includegraphics[height=6cm,width=8cm]{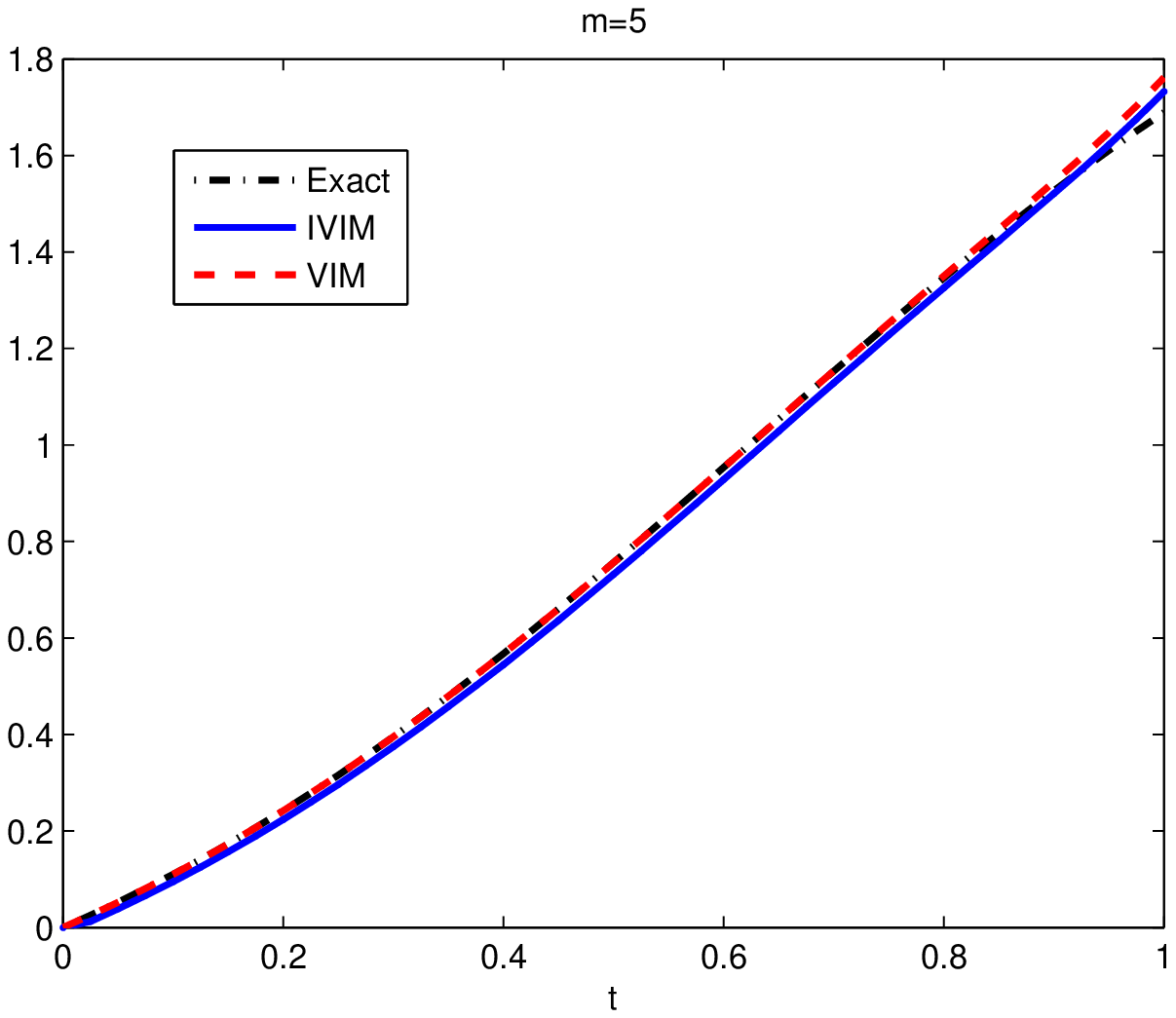}\includegraphics[height=6cm,width=8cm]{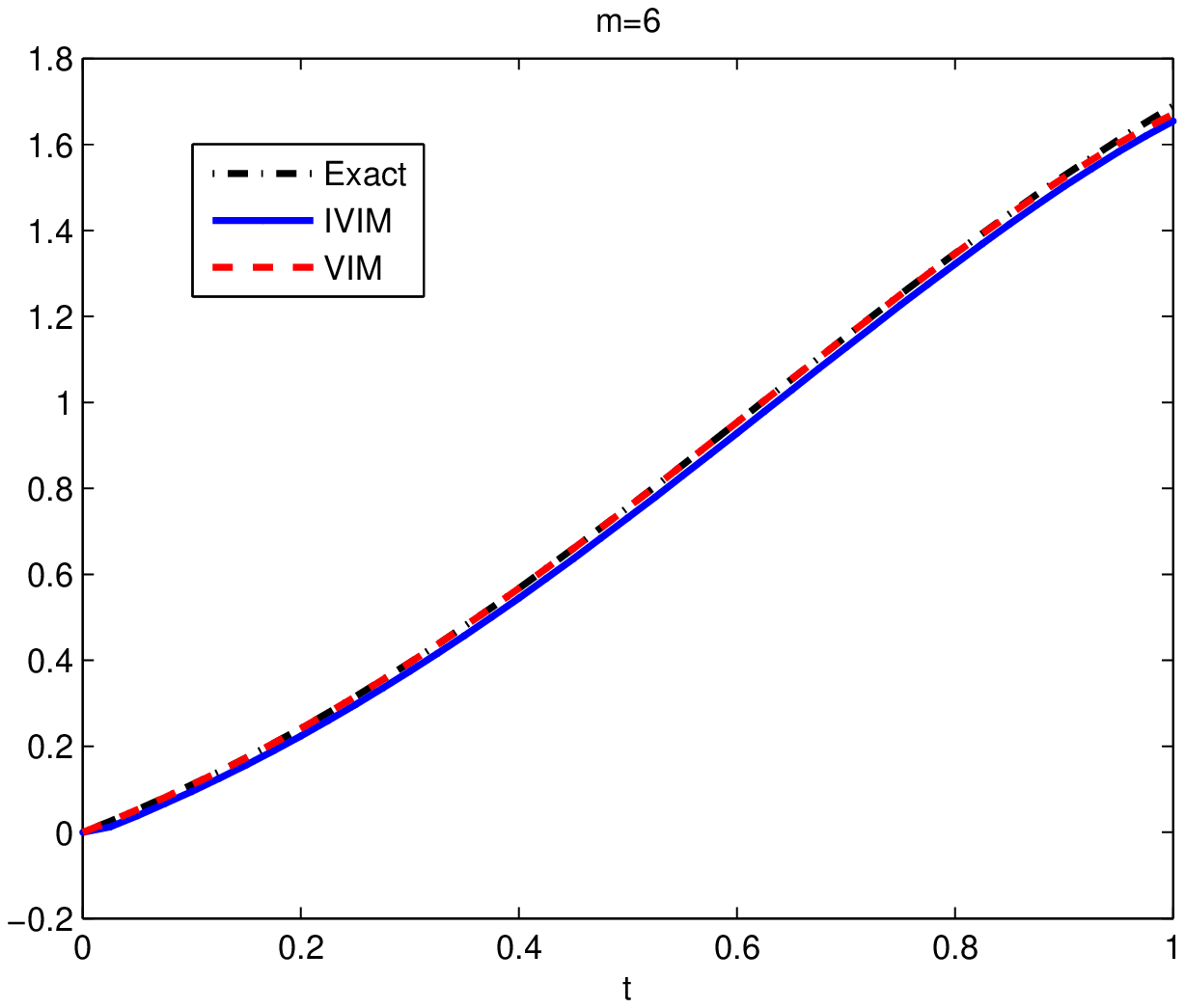}
\caption{A comparison between the exact and the approximate solutions computed by the VIM and the IVIM for the Riccati differential equation for Example 1.}\label{Fig1-Ex1}
\end{figure}

\noindent{\bf Example 1.} We consider the quadratic Riccati differential equation of the form
\begin{equation}\label{Eq012}
\left\{
  \begin{array}{ll}
    u'(t)=2u(t)-u^2(t)+1, & t\in [0,1], \\
    u(0)=0.
  \end{array}
\right.
\end{equation}
The exact solution of (\ref{Eq012}) is given by (see for example \cite{Abbasbandy})
\[
u(t)=1+\sqrt{2} \tanh \left(\sqrt{2}t+\frac{1}{2} \log \frac{\sqrt{2}-1}{\sqrt{2}+1} \right).
\]
We assume that $\mathcal{L}u=u'-2u$, $\mathcal{N}u= u^2$ and $g\equiv 1$. By setting $\alpha=-2$, according to Remark 1 we get $G_m\equiv 0$ and
\[
H_m(s,t)=e^{2(t-s)}(u^2_m(s)-1),
\]
for $(s,t)\in [0,1]\times [0,1]$. Therefore, Eq. (\ref{Eq005}) takes the following form
\[
u_{m+1}(t)=\int_{0}^t e^{2(t-s)}(1-u^2_m(s))ds,
\]
which is equivalent to the iteration produced by the VIM. In Figure \ref{Fig1-Ex1}, the approximate solutions computed  by the VIM and the IVIM  for different values of $m$ together with the exact solution are depicted. {In the IVIM},  we have considered $n=41$. The figure shows that the solution computed by the IVIM is in good agreement with that of the VIM and both of them converge to the exact solution. Table \ref{Tbl1} shows a comparison between the CPU times (in seconds) for computing the approximate solutions computed by these two methods. As we see, the CPU times for computing the approximated solution obtained by the VIM increase drastically as $m$ increases, whereas those of the IVIM are negligible.

\begin{figure}
\centering
\includegraphics[height=9cm,width=12cm]{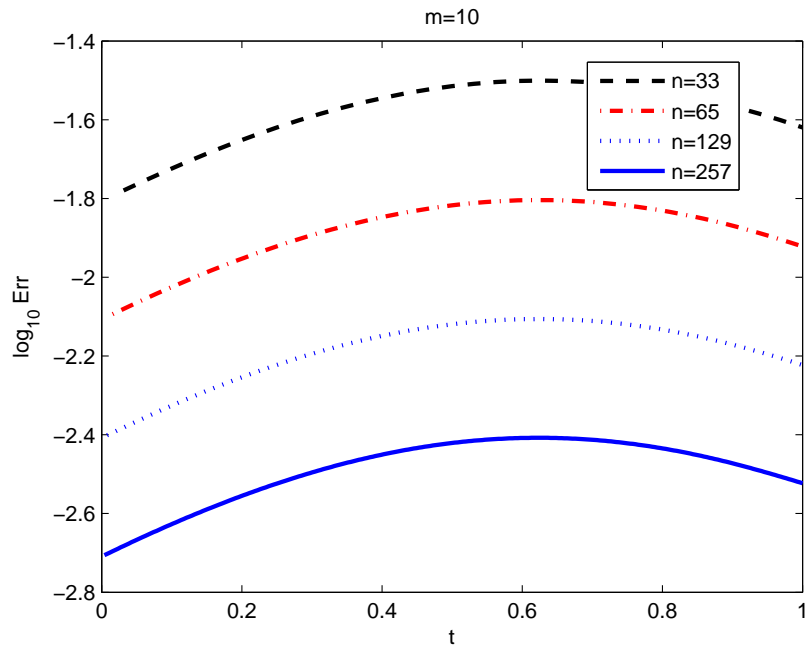}
\caption{$\log_{10}$ of the absolute error of the solutions provided by the IVIM for $m=10$ and $n=33,65,129,257$ in Example 1.}\label{Fig2-Ex1}
\end{figure}

\begin{figure}
\centering
\includegraphics[height=9cm,width=12cm]{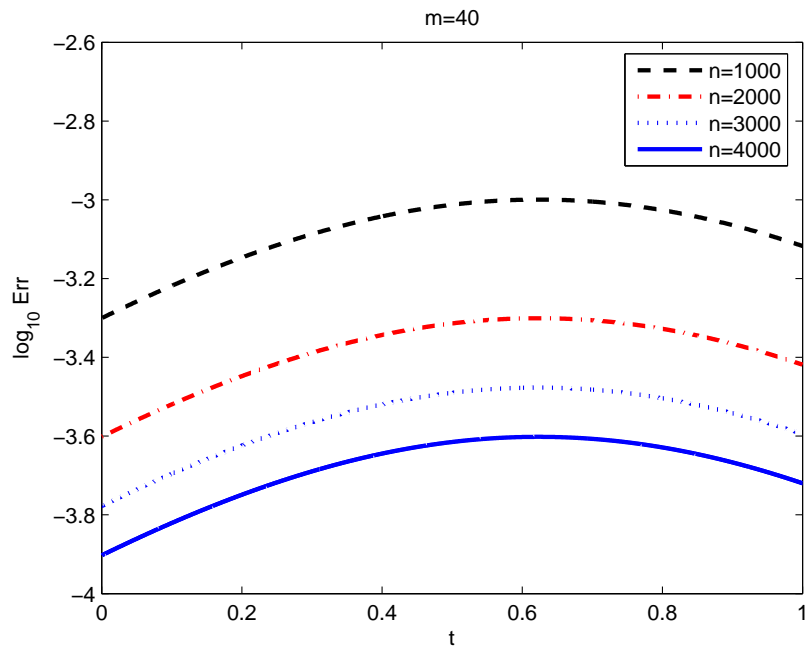}
\caption{$\log_{10}$ of the absolute error of the solutions provided by the IVIM for $m=10$ and $n=1000,2000,3000,4000$ in Example 1.}\label{Fig3-Ex1}
\end{figure}

For more investigation we set $m=10$ and $n=33,65,129,257$. In Figure \ref{Fig2-Ex1} the $\log_{10}$ of the absolute error of the computed solutions  by the IVIM,  denoted by $\log_{10} {\rm Err}$, are displayed. The CPU times for all of the four values of $n$ are $0.000$. Similarly, for  $m=40$ and $n=1000,2000,3000,4000$, the results are shown in Figure \ref{Fig3-Ex1}. Here the CPU times for computing the approximate solutions by the IVIM are $0.265,0.859,1.891,3.203$, respectively, for $n=1000,2000,3000,4000$. These figures show that the computed solutions obtained by the IVIM are in good agreement with the exact solution. Moreover, the results show that the CPU times for setting up the IVIM is very small in comparison to the VIM.

\begin{table}
\caption{CPU times (in seconds) for computing the approximate solutions by the VIM and the IVIM for Example 1.}
\vspace{-0.2cm}
\begin{center}
\begin{tabular}{lcccccccccc}\hline
$m$ \qquad\qquad  &  1       &   2       &    3        &     4     &  5       & 6      &  7      &  8 \\ \hline
VIM               &  0.015   &  0.093    &  0.109      &  0.204   &  0.672   & 2.485  &  20.034 &  383.469\\
IVIM              &  0.000   &  0.000    &  0.000      &  0.000   &  0.000   & 0.000  &  0.000  &  0.000\\ \hline
\end{tabular}
\end{center}
\label{Tbl1}
\end{table}

\bigskip

\begin{figure}
\centering
\includegraphics[height=9cm,width=12cm]{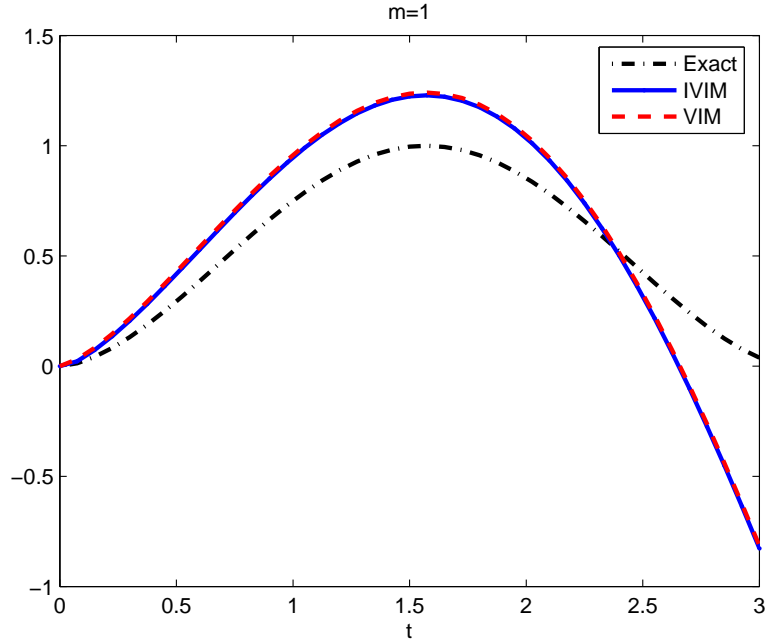}
\caption{Exact solution together with the approximate solutions computed by one iteration of the VIM and the IVIM for Example 2.}\label{Fig1-Ex2}
\end{figure}

\noindent{\bf Example 2.} We consider the initial value problem
\begin{equation}\label{Eq013}
\left\{
  \begin{array}{ll}
    u'(t)=\frac{5}{3} \sqrt[5]{u^2(t)}  \cos t, & t\in [0,3], \\
    u(0)=0.
  \end{array}
\right.
\end{equation}
The exact solution of the problem  is $u(t)= t\sqrt[3]{t^2} \sin t$. Let $\mathcal{L}u= u'$, $\mathcal{N}u(t)=-\frac{5}{3} \sqrt[5]{u^2(t)} \cos t$ and $g\equiv 0$. Setting $\alpha=0$, from Remark 1 we get $G_m\equiv 0$ and
\[
H_m(s,t)=-\frac{5}{3} \sqrt[5]{u_m^2(s)} \cos s,
\]
for $(s,t)\in [0,3]\times [0,3]$. Therefore, from Eq. (\ref{Eq005}) the VIM takes the following form
\[
u_{m+1}(t)=\int_{0}^t \frac{5}{3} \sqrt[5]{u_m^2(s)}   \cos s~ds.
\]
We have provided a MAPLE code for computing $u_m$'s. In Figure \ref{Fig1-Ex2} we have depicted the exact solution together with the computed solutions by one iteration of the VIM and IVIM. As is seen, there is a good agreement between the computed solutions  by the VIM and IVIM. For $m\geq 2$ the VIM could not provide an explicit expression for the $u_m$'s, while the IVIM is well suited  for sufficiently large values of $m$ and $n$ (Figures 5 and 6).

\begin{figure}
\centering
\includegraphics[height=9cm,width=12cm]{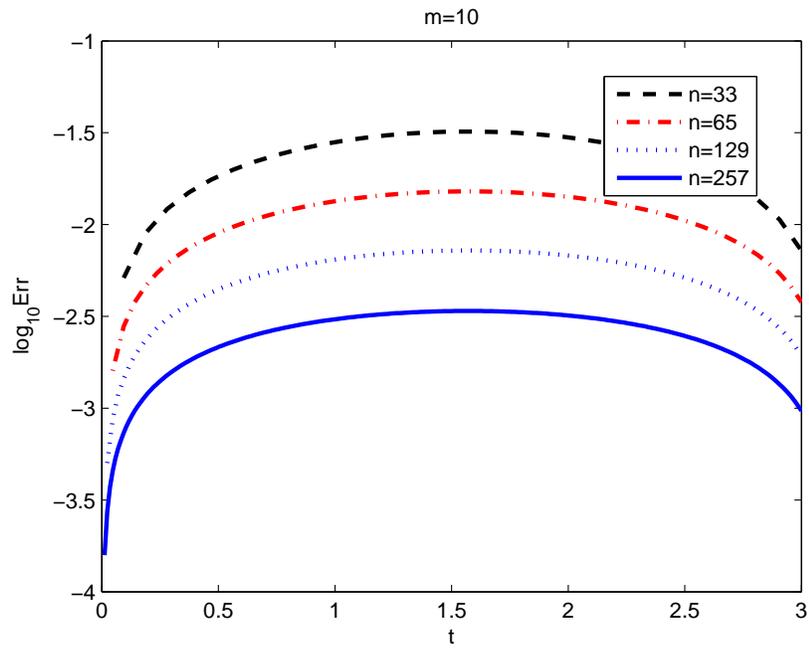}\caption{$\log_{10}$ of the absolute error of the computed solutions by the IVIM for $m=10$ and $n=33,65,129,257$ for Example 2.}\label{Fig2-Ex2}
\end{figure}

\begin{figure}
\centering
\includegraphics[height=9cm,width=12cm]{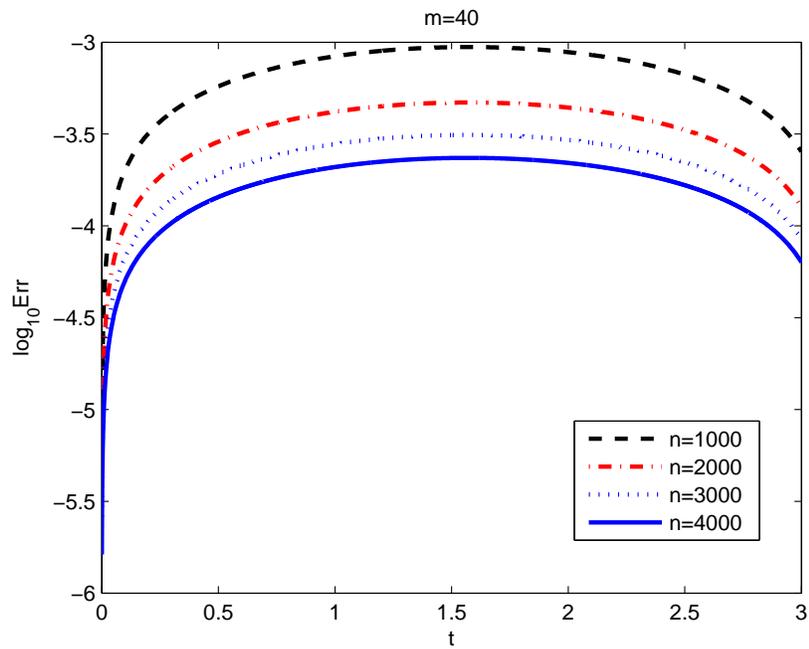}\caption{$\log_{10}$ of the absolute error of the computed solutions by the IVIM for $m=40$ and $n=1000,2000,3000,4000$ for Example 2.}\label{Fig3-Ex2}
\end{figure}

To see the efficiency of the IVIM, we consider the solutions provided by the IVIM with $m=10$ and different values of $n$ ($=33,65,129$ and $257$).  Figure \ref{Fig2-Ex2} shows the $\log_{10}$ of the absolute error for the given solutions by IVIM. As we can see, the CPU times for computing these solutions are $0.000, 0.016,0.031$ and $0.047$, respectively. Similar results have been shown in Figure \ref{Fig3-Ex2} for $m=40$ and  $n=1000,2000,3000$ and $4000$. The CPU times for computing the approximate solutions are $0.281,0.891,2.234$ and $3.484$, respectively.

\bigskip

\noindent{\bf Example 3.} Although we have presented the IVIM for the first-order initial value problems, by this  example we illustrate that the method can be implemented for the higher-order initial value problems. To do so, consider the second-order initial value problem
\begin{equation}\label{Ex3Eq1}
\left\{
  \begin{array}{ll}
    u''(t)-2u'(t)^2+u'(t)+u(t)=g(t), & t\in [0,1.5], \\
    u(0)=u'(0)=0, &
  \end{array}
\right.
\end{equation}
where $g(t)=t+2\sin^2 \frac{t}{2}-8\sin^4 \frac{t}{2}$. The exact  solution of this problem is $u(t)=x-\sin t$. Letting $v(t)=u'(t)$ for $t\in [0,1.5]$, Eq. (\ref{Ex3Eq1}) can be written as a
system of two simultaneous first-order differential equations
\begin{equation}\label{Ex3Eq2}
\left\{
  \begin{array}{ll}
    u'(t)=v(t),\\
    v'(t)-2v(t)^2+v(t)+u(t)=g(t),
  \end{array}
\right.
\end{equation}
with the initial condition $u(0)=v(0)=0$. Having in mind the exact  solution of (\ref{Ex3Eq1}), we see that the exact solution of the problem (\ref{Ex3Eq2}) is given by $(u(t),v(t))=(t-\sin t, 1-\cos t)$.

To implement the VIM for solving (\ref{Ex3Eq2}),  in the first equation of (\ref{Ex3Eq2}) the linear and nonlinear terms are chosen as $u'$ and $-v$, respectively. In the same way, for the second equation in  (\ref{Ex3Eq2}),  $v'+v$ and $-2v^2+u$ are considered as the linear and nonlinear terms, respectively.
  Similar to the previous examples the VIM iteration takes the following form (for more details the reader can refer to \cite{Zhao}):
\begin{equation}\label{Ex3Eq3}
\left\{
  \begin{array}{ll}
    u_{m+1}(t)=u_m(t)-\displaystyle\int_{0}^{t} \left( u_m'(s)-v_m(s) \right) ds,\\[4mm]
    v_{m+1}(t)=v_m(t)-\displaystyle\int_{0}^{t} e^{s-t}\left( v_m'(s)-2v_m(s)^2+v_m(s)+u_m(s)-g(s)\right)ds,
  \end{array}
\right.
\end{equation}
where $u_0$ and $v_0$ are two given functions satisfying $u_0(0)=v_0(0)=0$. Obviously, from (\ref{Ex3Eq3}) it follows that  $u_m(0)=v_m(0)=0$. Therefore, by straightforward integration by parts  Eq. (\ref{Ex3Eq3}) can be written in the simpler form
\begin{equation}\label{Ex3Eq4}
\left\{
  \begin{array}{ll}
    u_{m+1}(t)=\displaystyle\int_{0}^{t} v_m(s)ds,\\[4mm]
    v_{m+1}(t)=\displaystyle\int_{0}^{t} e^{s-t}\left(2v_m(s)^2-u_m(s)+g(s)\right)ds\\[4mm]
             \hspace{1.45cm}=\displaystyle\int_{0}^{t} e^{s-t}\left(2v_m(s)^2+g(s)\right)ds
             -\displaystyle\int_{0}^{t} e^{s-t}u_m(s)ds .
  \end{array}
\right.
\end{equation}
To implement the IVIM, we set
\begin{eqnarray*}
\bar{H}_m(s,t) &=& -v_m(s), \\
\hat{H}_m(s,t) &=& -e^{s-t}\left(2v_m(s)^2+g(s)\right), \\
\tilde{H}_m(s,t) &=& e^{s-t}u_m(s).
\end{eqnarray*}
In this case, Eq.  (\ref{Ex3Eq4}) can be rewritten as
\begin{equation}\label{Ex3Eq5}
\left\{
  \begin{array}{ll}
    u_{m+1}(t)=-\displaystyle\int_{0}^{t} \bar{H}_m(s,t)ds,\\[4mm]
    v_{m+1}(t)=-\displaystyle\int_{0}^{t} \hat{H}_m(s,t) ds-\displaystyle\int_{0}^{t} \tilde{H}_m(s,t) ds .
  \end{array}
\right.
\end{equation}
Taking Eq. (\ref{Eq005}) into account, one can approximate $u_{m+1}$ and $v_{m+1}$ by the idea described in Section \ref{SEC3}. Since the implementation of the method is straightforward, the details of the method are omitted here. Since, $\hat{H}_m$ involves the nonlinear term $2v_m^2$, the main difficulties mentioned in Example 1 arise here, too.

\begin{figure}
\centering
\includegraphics[height=4.5cm,width=7cm]{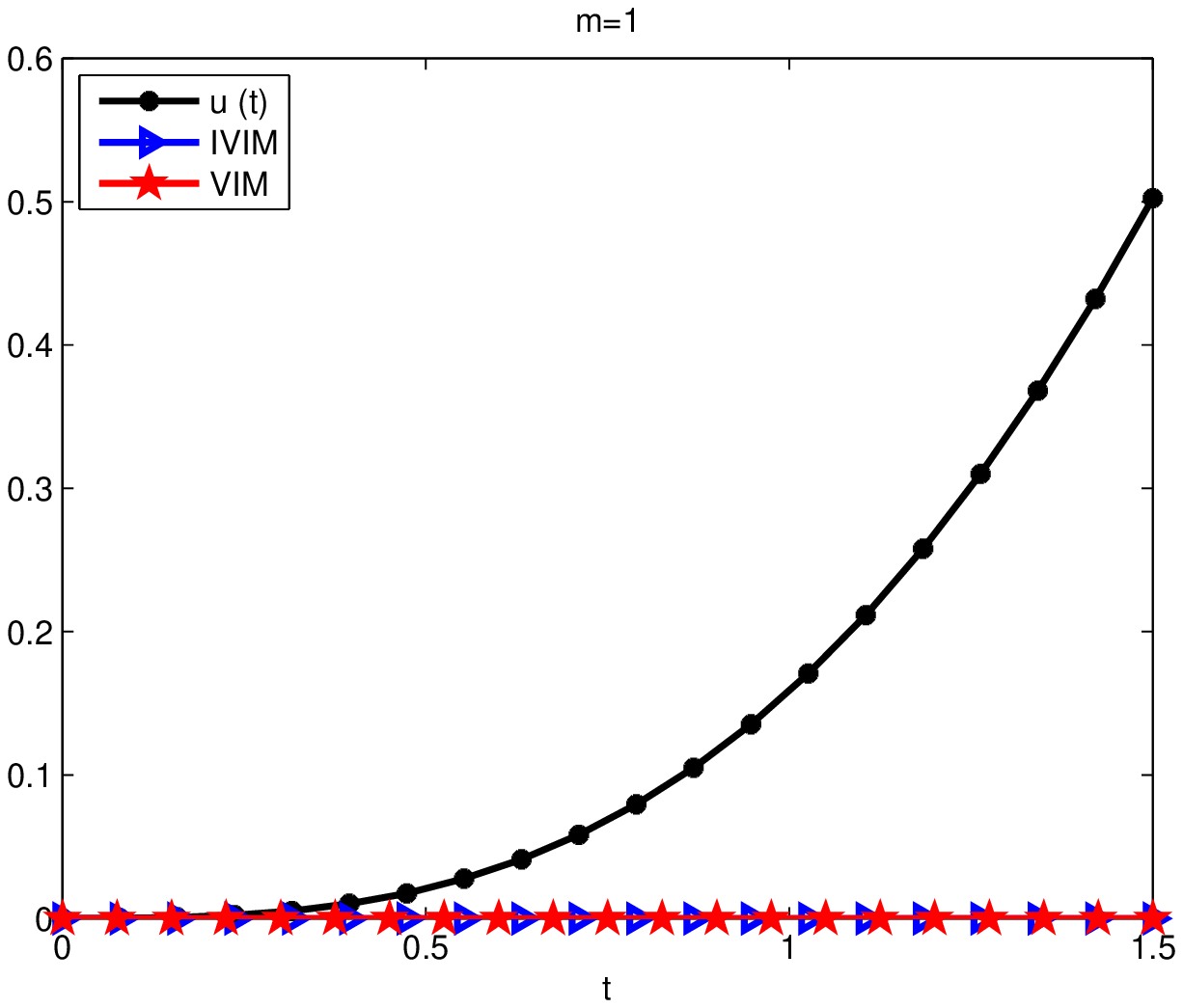}\includegraphics[height=4.5cm,width=7cm]{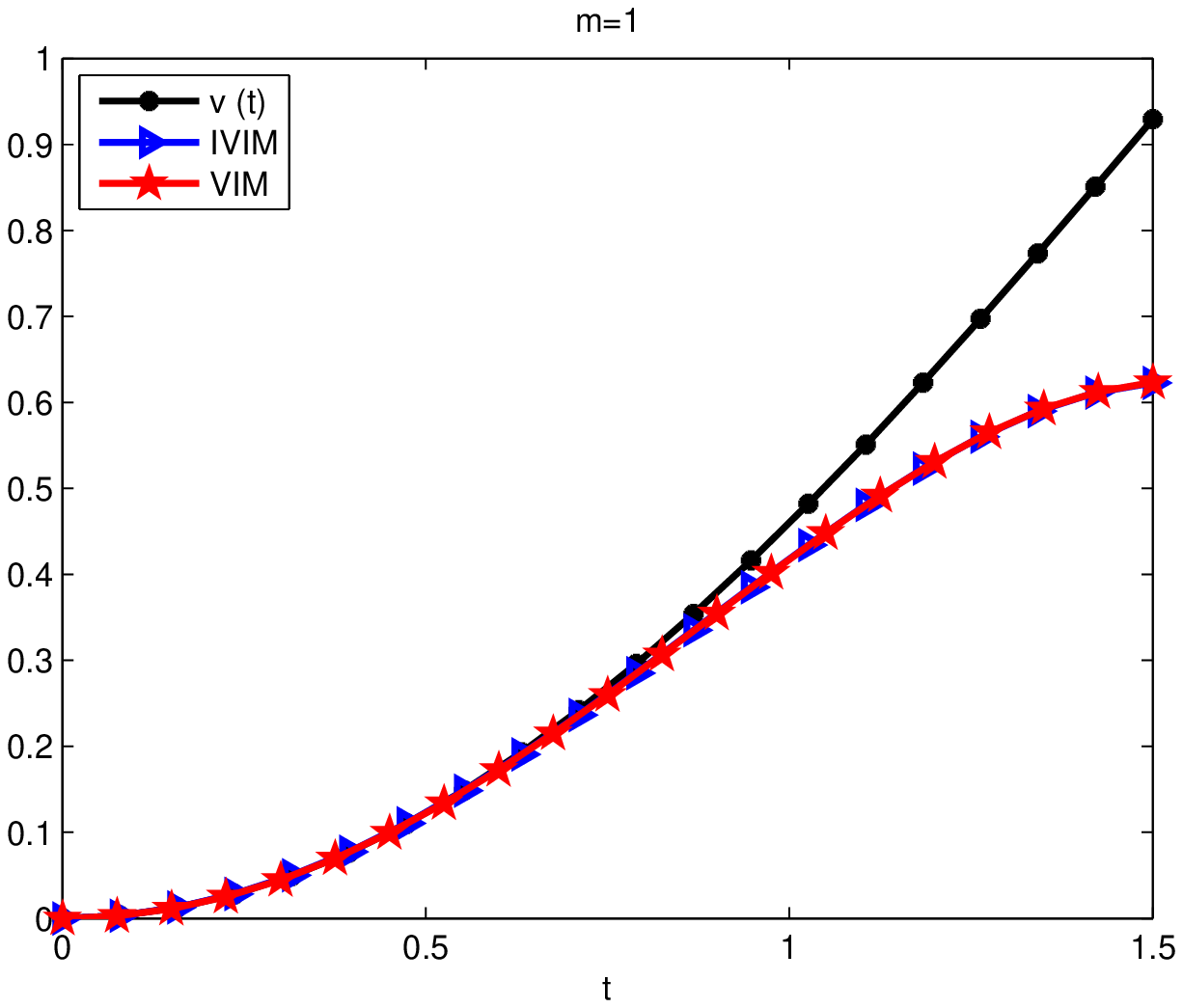}
\includegraphics[height=4.5cm,width=7cm]{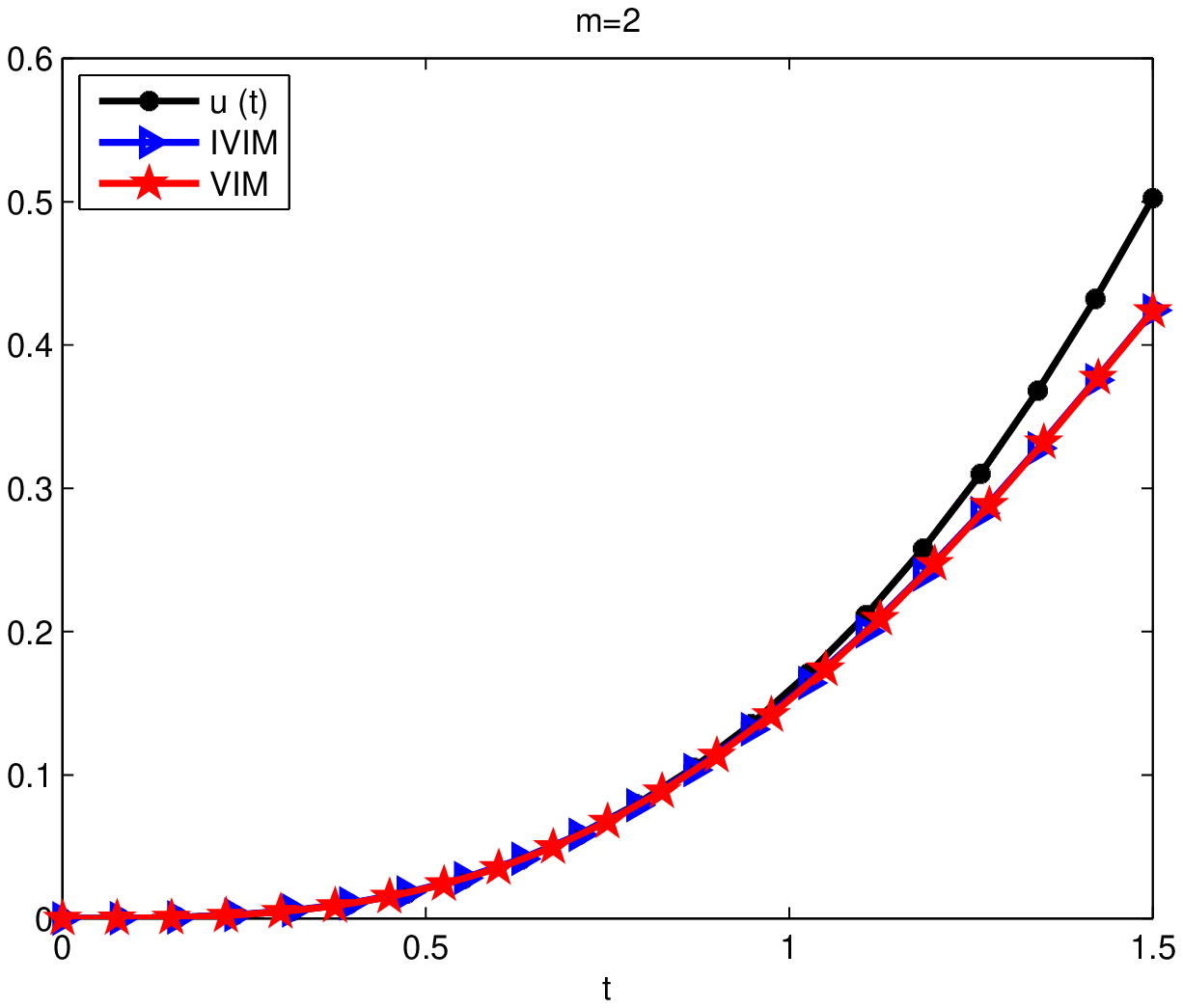}\includegraphics[height=4.5cm,width=7cm]{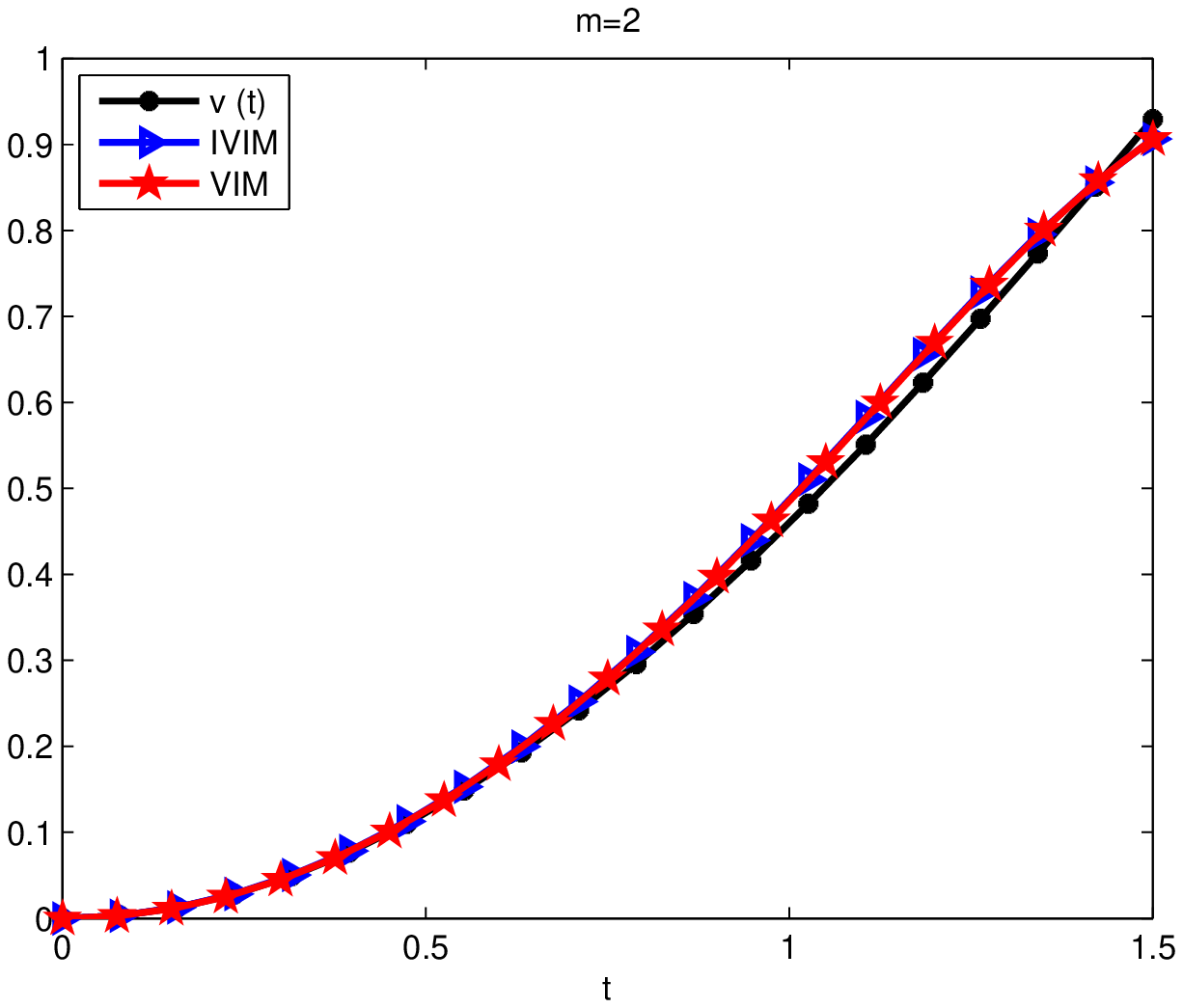}
\includegraphics[height=4.5cm,width=7cm]{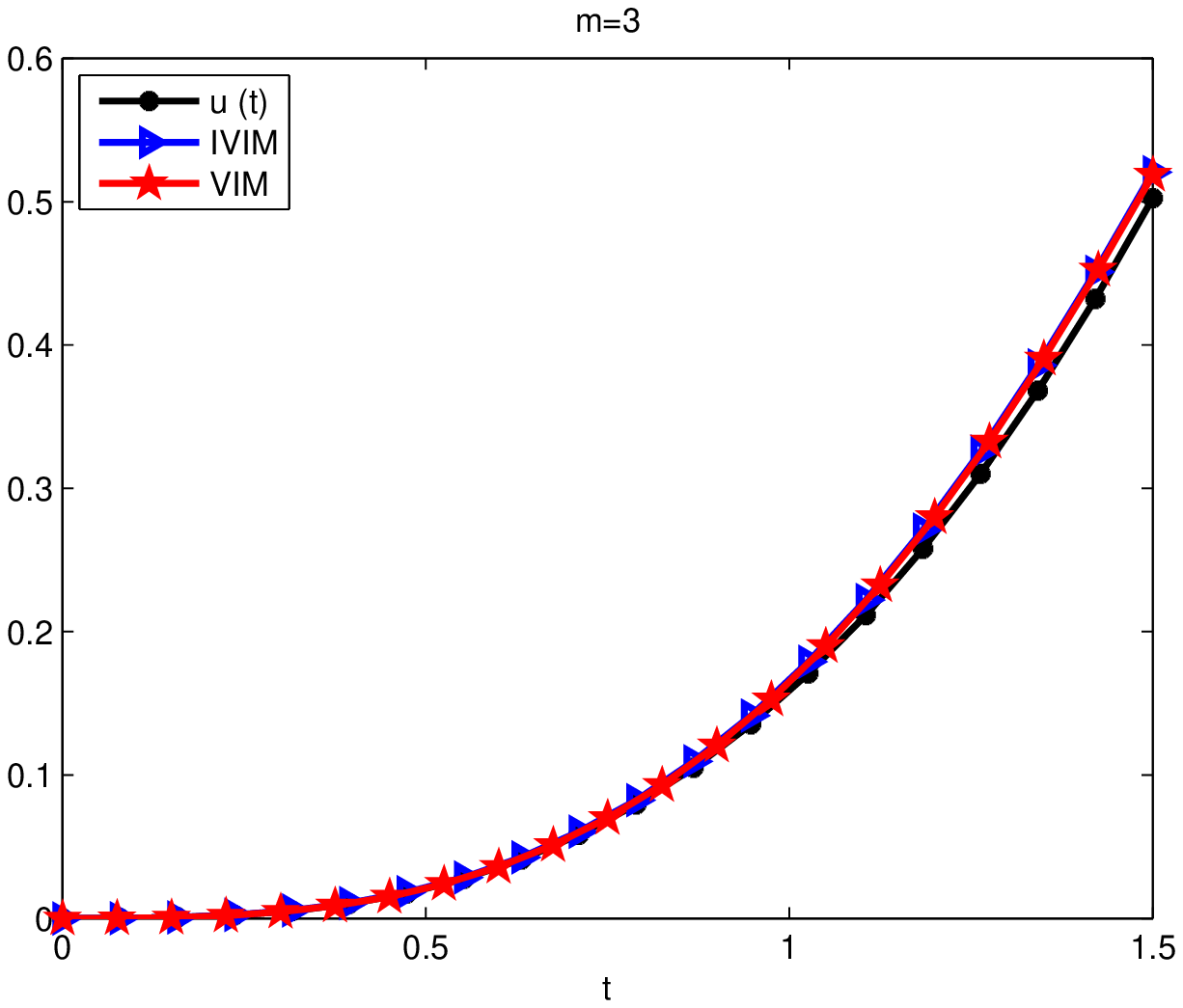}\includegraphics[height=4.5cm,width=7cm]{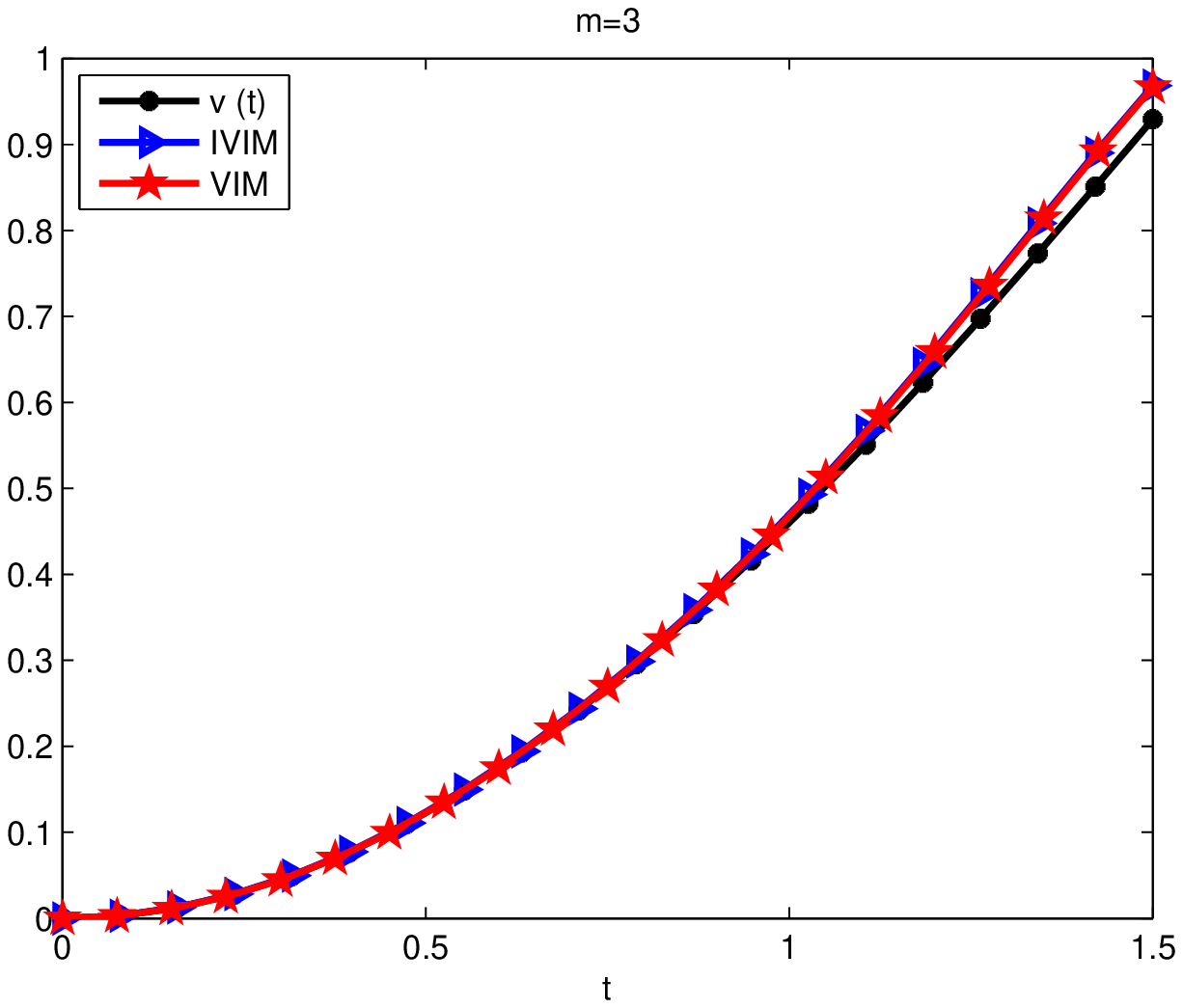}
\includegraphics[height=4.5cm,width=7cm]{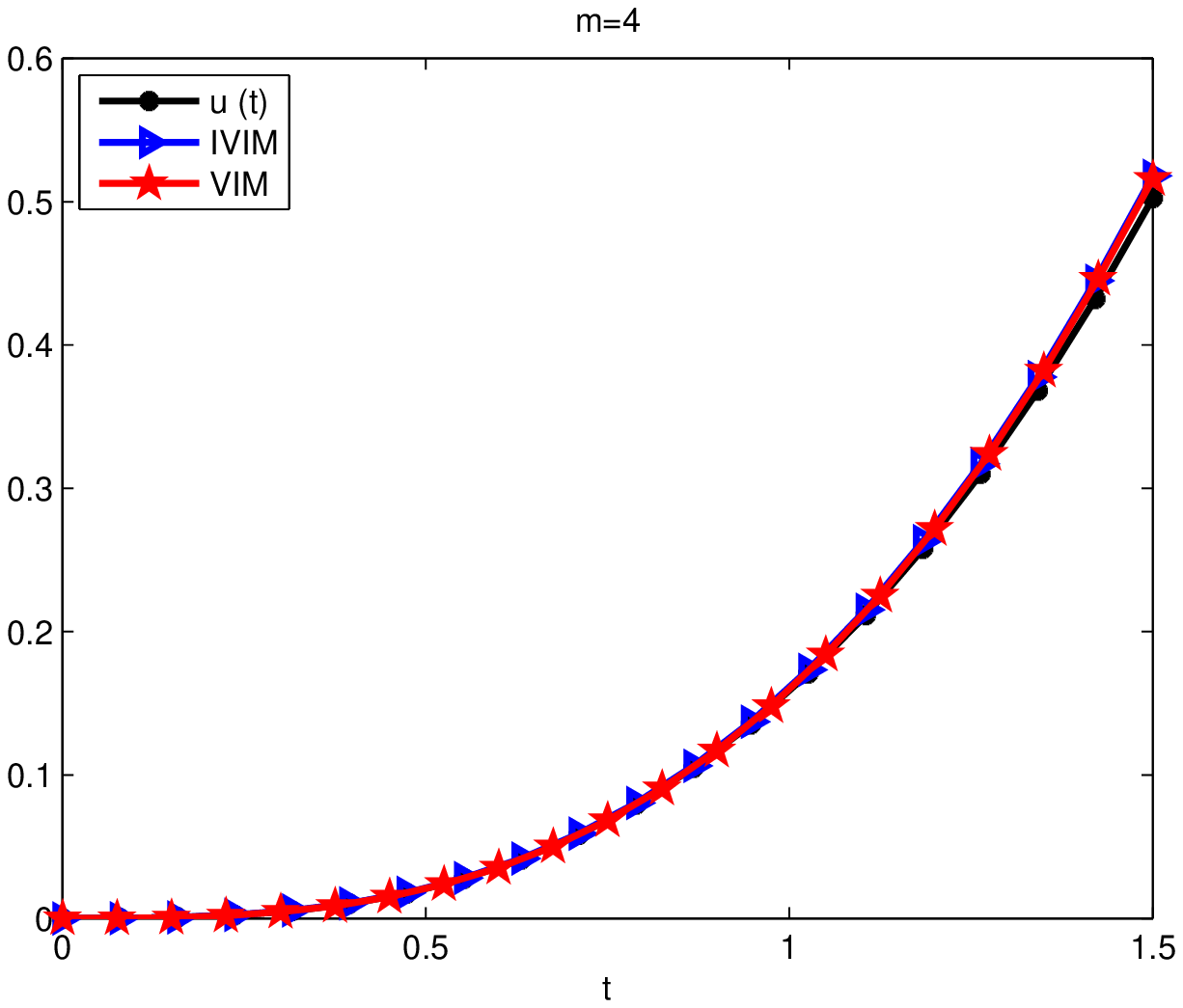}\includegraphics[height=4.5cm,width=7cm]{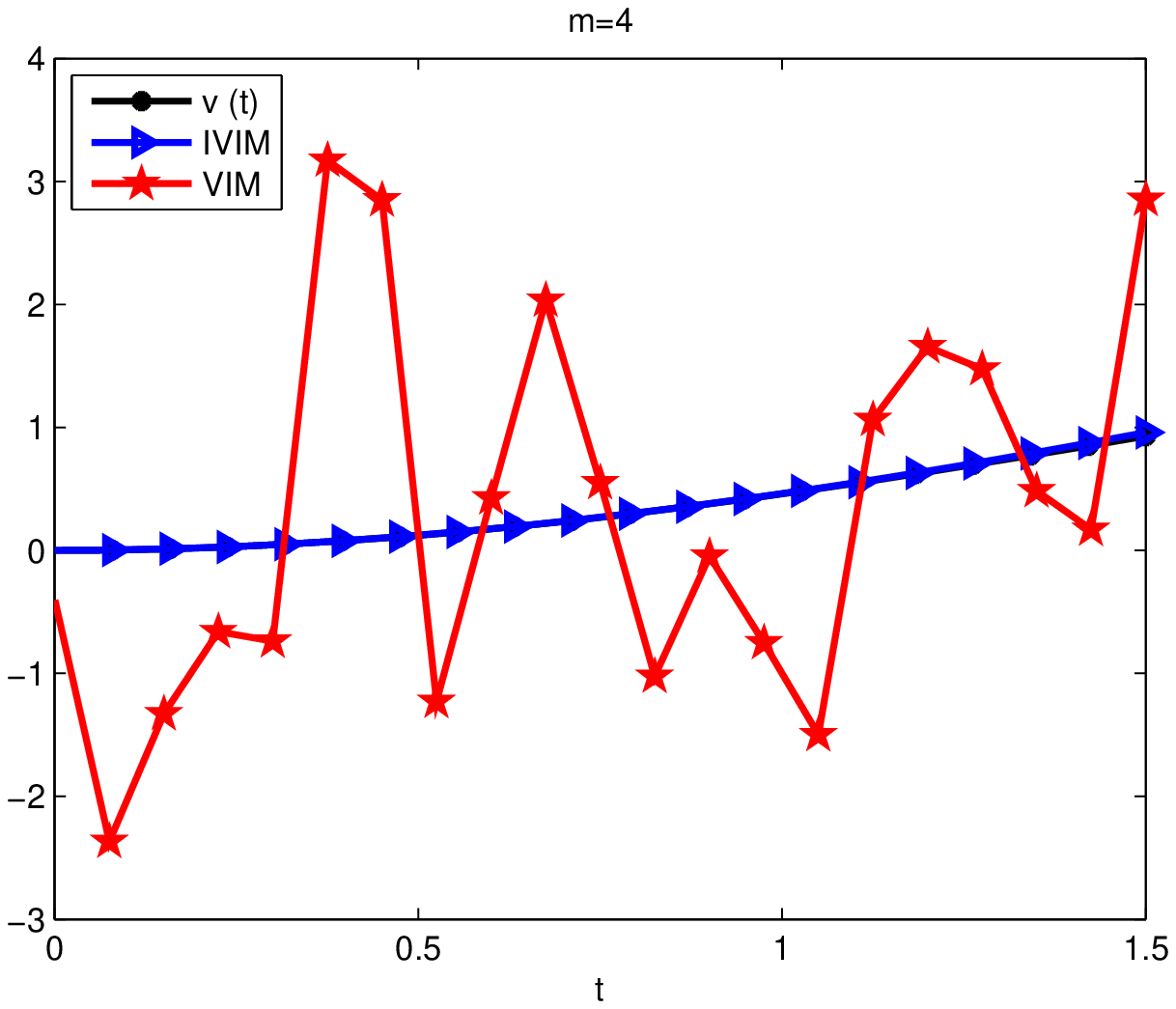}
\includegraphics[height=4.5cm,width=7cm]{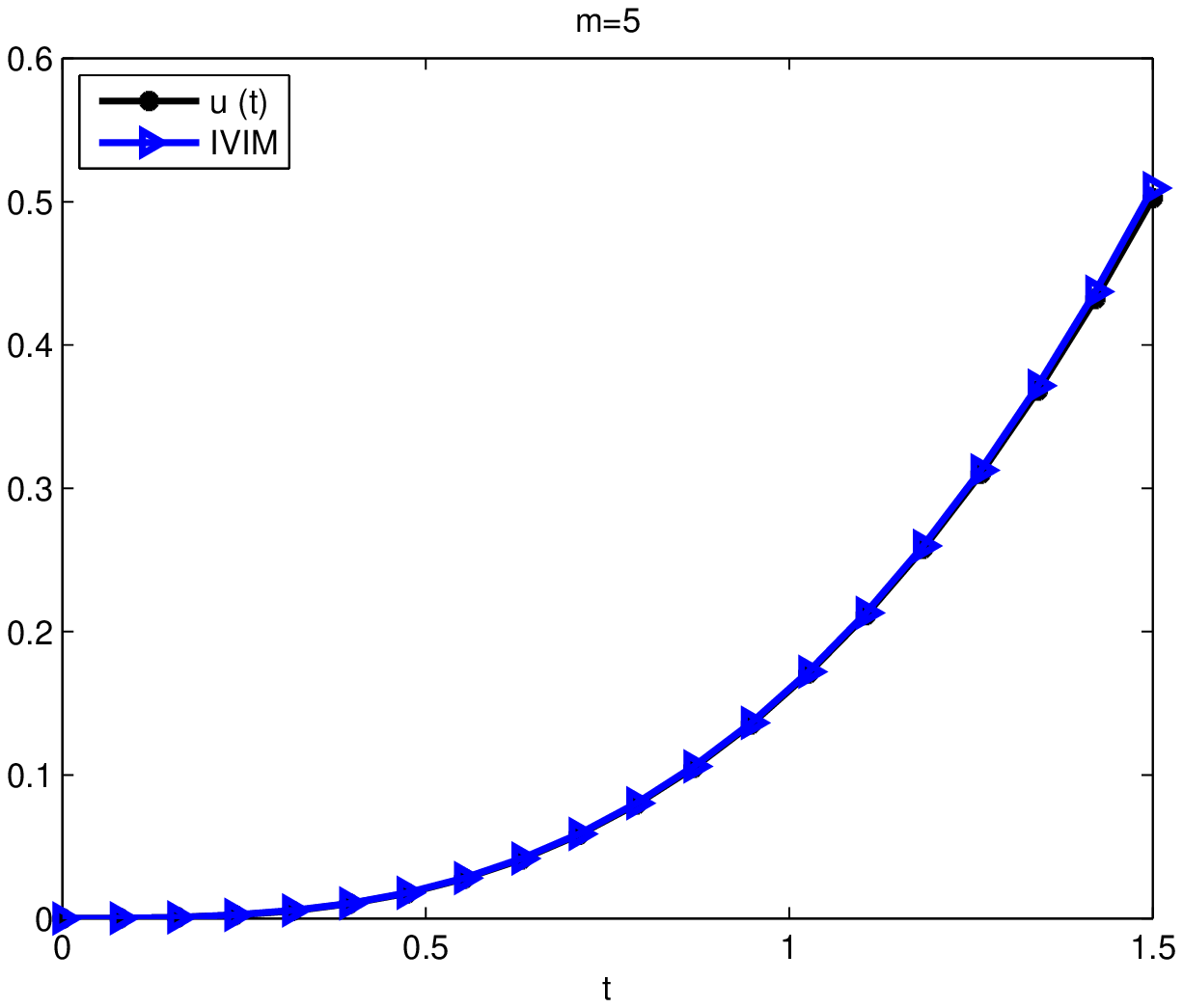}\includegraphics[height=4.5cm,width=7cm]{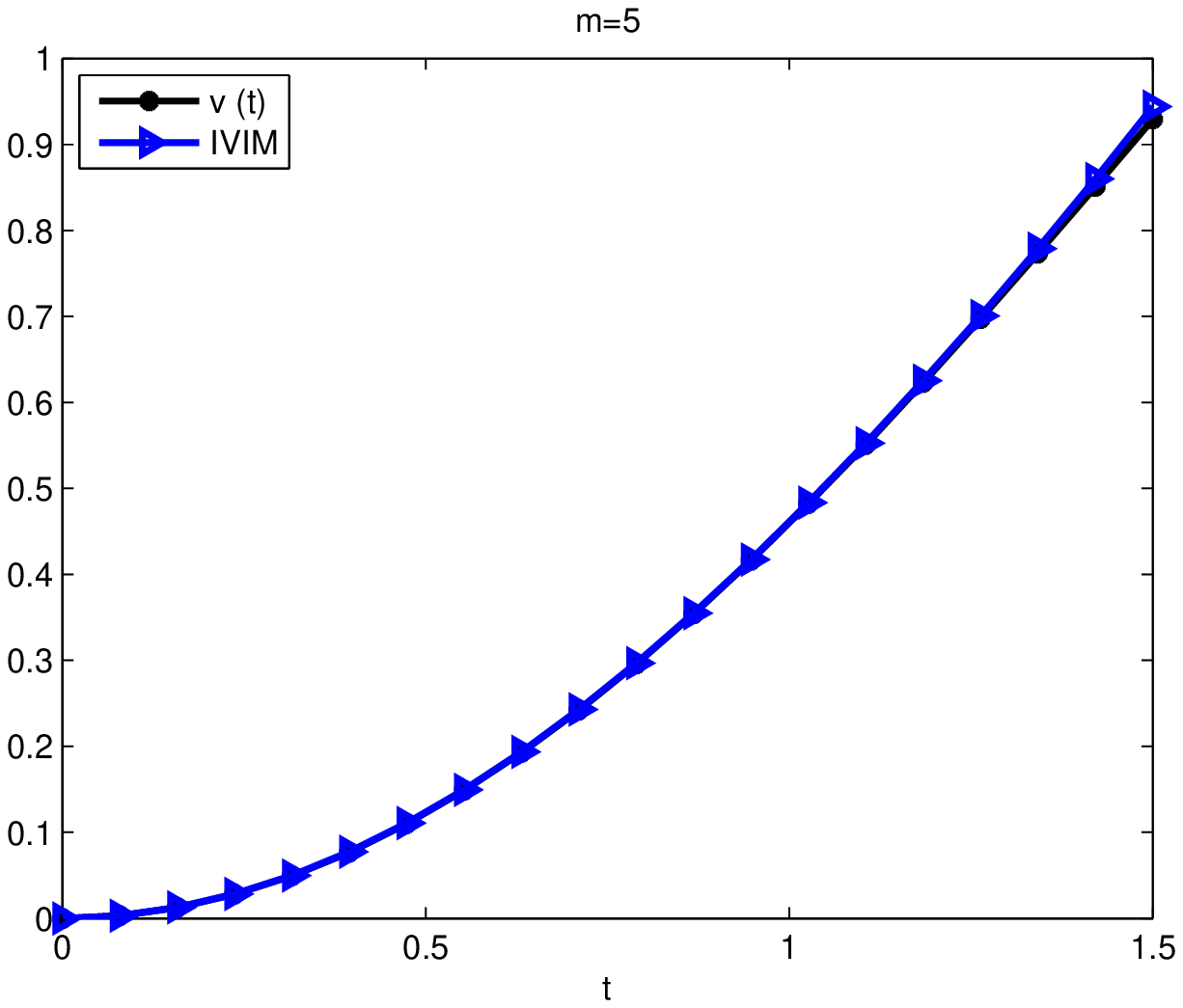}
\caption{A comparison between the graphs of $u(t)=t-\sin t$ and $v(t)=u'(t)=1-\cos(t)$ , respectively, with $u_m(t)$ and $v_m(t)$, for $m=1,2,3,4,5$.}\label{Ex3Fig1}
\end{figure}

For the numerical results we use $u_0(t)=v_0(t)=0$ for $t\in [0,1.5]$ as the starting point. We now present the numerical results of using the IVIM for solving the problem. To this end, we first take $n=20$. In Figure \ref{Ex3Fig1} the graphs of $u_{m}$ and $v_{m}$ for $m=1,2,3,4,5$, computed by the VIM and the IVIM, together with $u(t)=t-\sin t$ and $v(t)=1-\cos t$ are depicted.  This figure illustrates that, in the IVIM,  $u_m(t)\rightarrow u(t)$ and $v_m(t)\rightarrow v(t)$  for $t\in [0,1.5]$ as $m$ tends to infinity. Moreover, it shows that there is a good agreement between the solution provided by the IVIM and the VIM when the VIM gives a suitable solution.   In addition, the VIM can not provides any suitable expression for  $v_4$, $u_m$ and $v_m$ for $m\geq 5$. This is due to the nonlinear term involving $v_{m+1}$. For more investigation we report the CPU times (in seconds) for computing the solutions by the VIM and the IVIM methods in Table \ref{Tbl2}. As we can observe, the CPU times for the VIM are very small, whereas those of the VIM drastically grows as $m$ increases.

\begin{table}
\caption{CPU times (in seconds) for computing the approximate solutions by the VIM and the IVIM for Example 3.}
\vspace{-0.2cm}
\begin{center}
\begin{tabular}{lcccccccccc}\hline
$m$ \qquad\qquad  &  1       &   2       &    3        &     4     &   5     \\ \hline
VIM               &  0.046   &  0.313    &  1.531      &  16.750   &  576.4 (Fail)  \\
IVIM              &  0.010   &  0.012    &  0.012      &  0.014    &  0.015  \\ \hline
\end{tabular}
\end{center}
\label{Tbl2}
\end{table}

\begin{figure}
\centering
\includegraphics[height=6.cm,width=8cm]{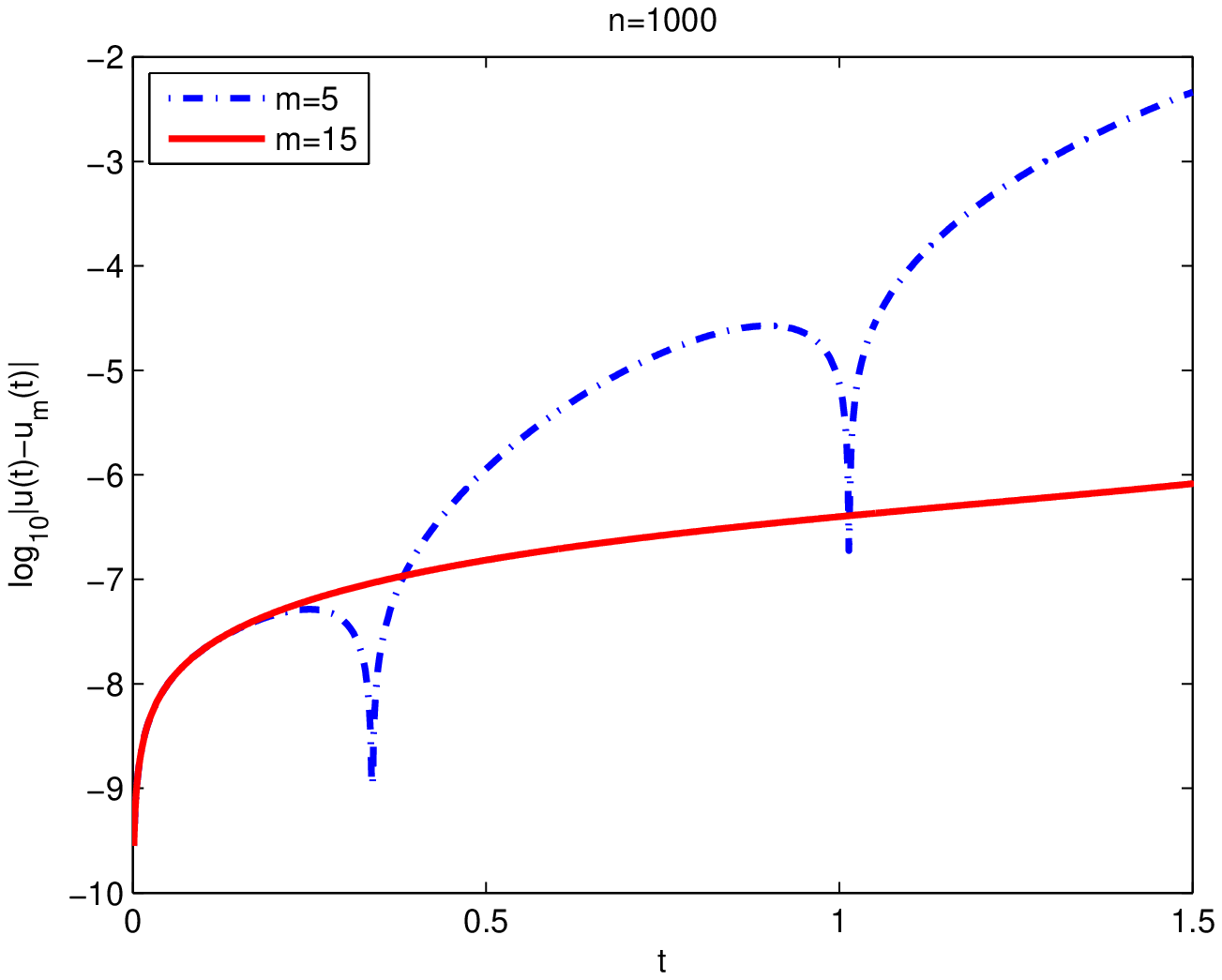}\includegraphics[height=6.cm,width=8cm]{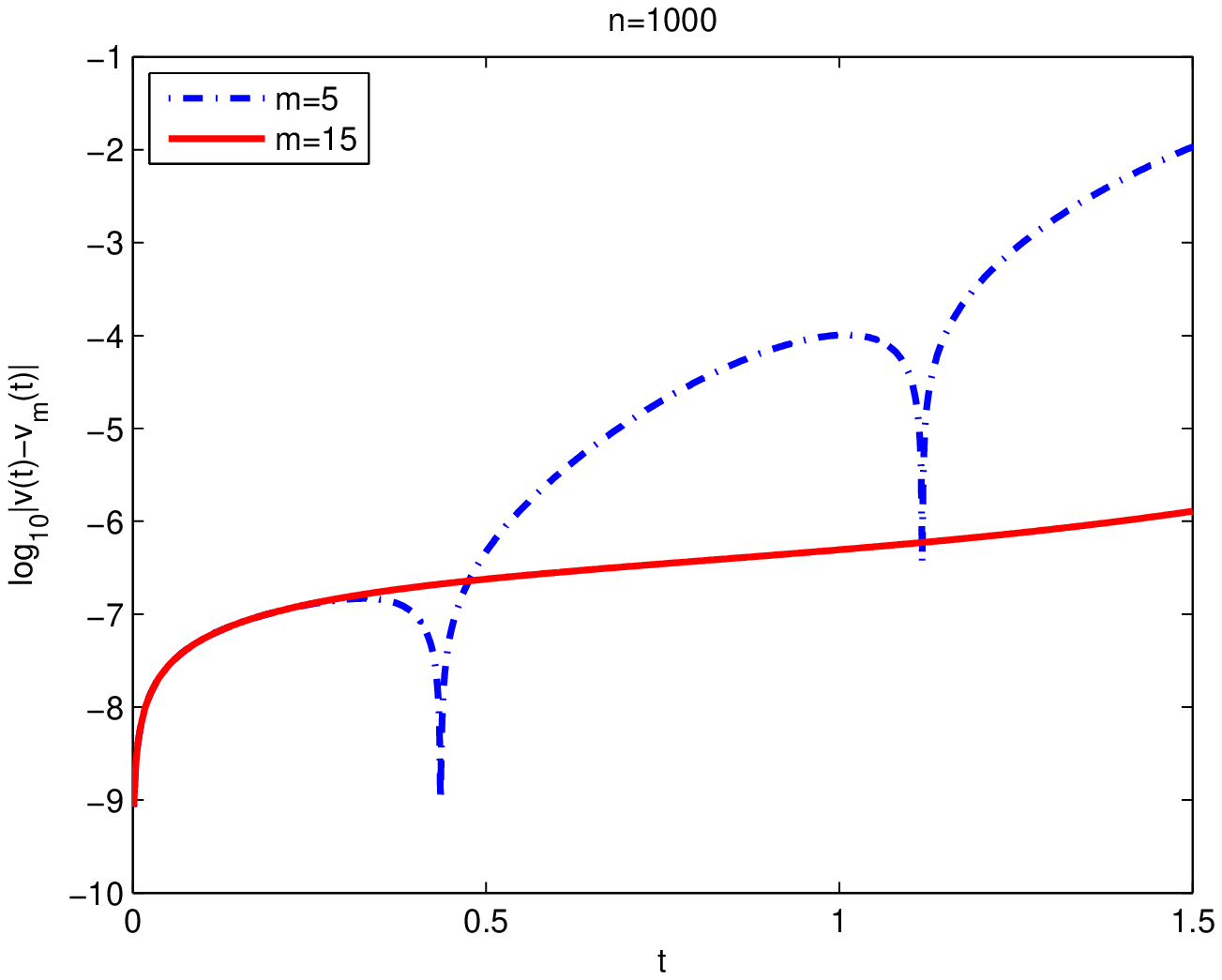}
\caption{   $\log_{10}$ of the absolute error of the solution computed by the IVIM for Example 3 for $n=1000$ and $m=5,15$.}\label{Fig2Ex3}
\end{figure}

Finally, we set $n=1000$ and compute the approximate solutions of the problem by the IVIM for $m=5$ and $m=15$. In Figure \ref{Fig2Ex3},  $\log_{10}$ of the absolute error of the computed solutions is displayed.  It is noted that the VIM gives the approximate solutions in 0.082 and 0.201 seconds for $m=5$ and $m=15$, respectively.  As we see, the IVIM  provides quite suitable solutions for the problem in a reasonable amount of time.

\section{Conclusion}\label{SEC5}

We have proposed the interpolated VIM (IVIM) for solving one-dimensional initial value problems. The convergence of the method together with some numerical examples have been investigated. The numerical results show that the IVIM is more effective than the VIM. We have also shown the method can be easily implemented to the systems of ODEs.
The following aspects can be considered in future works:
\begin{enumerate}
\item Generalization of the method to two- or three-dimensional initial value problems;
\item Implementation of the idea  for the other methods such as the homotopy perturbation method;
\item  Using another type of basis functions such as the cubic spline functions and the Chebyshev polynomials.
\end{enumerate}

\section*{Acknowledgments}

The authors would like to thank the anonymous referees and the editors of the journal for their valuable comments and suggestions which substantially improved the quality of the paper. The work of the first author is partially supported by University of Guilan.



\end{document}